\documentclass[11pt]{article}
\usepackage[margin=3.1cm,top=3.3cm,bottom=3.1cm]{geometry}
\usepackage{amssymb, amsmath, amsthm, amsfonts, amscd}
\usepackage[english]{babel}
\usepackage{graphicx}
\usepackage{caption}
\usepackage{subcaption}
\usepackage{color,soul}
\usepackage{booktabs,dcolumn}                     
\usepackage[figureposition=bottom]{caption}       
\usepackage[usenames,dvipsnames,svgnames,table]{xcolor}
\usepackage[section,subsection,subsubsection]{extraplaceins}

\newcolumntype{d}[1]{D{.}{.}{#1}}

\newtheorem{theorem}{Theorem}[section]

\newtheorem{definition}[theorem]{Definition}

\numberwithin{equation}{section}
\newtheorem{remark}{Remark}

\def\N{\mathbb{N}}

\definecolor{orange}{rgb}{1,0.3,0.6}

\newcommand{\set}[1]{\{ #1 \}}

\title{Approximations of steady periodic water waves in flows with constant vorticity.}
\author{A.~Constantin\thanks{Faculty of Mathematics, University of Vienna, Oskar-Morgenstern-Platz 1, 1090 Vienna, Austria} \and
K.~Kalimeris\thanks{Radon Institute of Computational and Applied Mathematics, Altenberger Str.~69, 4040 Linz, Austria} \and 
O.~Scherzer \thanks{Computational Science Center, University of Vienna, Oskar-Morgenstern-Platz~1, 1090 Vienna, Austria and 
Radon Institute of Computational and Applied Mathematics, Altenberger Str.~69, 4040 Linz, Austria}
}

\begin{document}

\maketitle

\begin{abstract}
We provide high-order approximations to periodic travelling wave profiles and to 
the velocity field and the pressure beneath the waves, in flows with constant vorticity over a flat bed.
\end{abstract}



\section{Introduction}

Winds in offshore storms transmit energy to the ocean surface, creating waves. 
Once the waves leave the storm area they become organised into two-dimensional regular wave trains, in 
the form of a regular profile that propagates practically unchanged in a fixed direction. These periodic 
travelling waves are called `swell' in oceanography. Their shape is influenced by the underlying currents. 
The most significant currents on areas of the continental shelf and in many coastal inlets are the tidal 
currents \cite{Jon}. These alternating horizontal movements of water are created by the gravitational 
pull of the moon, and to a lesser degree, the sun, on the earth's surface, being 
associated with the rise and fall of the tide: the current associated with a rising 
water level, called the flood, is directed towards the shore, while the current associated with
a receding water level, called the ebb, is directed back out to sea. Flows of constant vorticity with a 
flat free surface provide adequate descriptions of pure tidal currents, cf. the discussion in 
\cite{Cb}, positive constant vorticity $\gamma>0$ being appropriate 
for the modelling of the ebb current and negative constant vorticity $\gamma<0$ for the flood current. 
A spectacular example of wave-current interaction is the Columbia River entrance, made by 
appreciable tidal currents one of the most hazardous navigational regions in the world since 
wave heights can easily be doubled in just a few hours \cite{Jon}. At this location, tidal velocities 
in excess of 2 m/s are encountered and in winter wave heights in excess 
of 6 m, up to 14-15 m, are common, cf. \cite{KO}. Note that moderate 
tidal currents reach speeds of up to 0.7 m/s, cf. \cite{WP}, while tidal currents with speeds 
of 5.5 m/s are encountered in between the Scottish mainland and the Orkney Islands, 
in the Pentland Firth --- see the data provided in \cite{G}.

In this paper we provide accurate approximations for 
the interaction of waves with underlying currents of constant vorticity 
in the absence of flow-reversal; for theoretical studies and numerical simulations 
of wave-current interactions with flow-reversal we refer to \cite{CSV, CV, W} and \cite{OS, V}, 
respectively. Let us specify that in a flow where both waves and currents are
present,  given velocity measurements at one point, the `current' is 
defined as the average velocity, and the periodic components that vary 
around this average are ascribed to the wave motion. Since irrotational flows 
do not present flow-reversals as the underlying current must 
be uniform, one of our main purposes will be to highlight the effects of vorticity on the surface wave profile. 
This problem is of great practical relevance since the detection of non-uniform underlying currents from 
the surface wave pattern is particularly important to anticipate and avoid hazardous conditions, 
where possible. In addition to the wave profile, we will also provide approximations for 
the velocity profile and 
the pressure throughout the fluid. Not only is there a strong interplay of these flow characteristics that 
is very useful in qualitative studies --- see \cite{C, CSp}, but in practice information on 
the state of the sea surface is often gathered from subsurface pressure and/or velocity 
measurements, cf. the discussions in \cite{Cl2, ClC, C2, K, OVD, VO}. Let us point out the following counter-intuitive 
fact: periodic travelling waves that propagate at the surface of water with a flat bed in a flow of constant 
vorticity must be symmetric if no flow-reversal occurs and if the wave profile is 
monotone between successive crests and troughs, cf. 
\cite{CEW, CE0, MM}. This means that an underlying non-uniform current of constant vorticity does 
not break the symmetry of irrotational wave trains, so that the manifestation of vorticity on the surface 
wave pattern must take subtler forms. 

The fact that even when both the wave motion and the underlying current of constant 
vorticity are known with accuracy, their interaction produces a significantly
different effect from that obtained by simply adding the effect of
the waves and the currents considered separately, cf. the discussions in \cite{SCJ, T, TK}, 
shows the importance of dealing simultaneously with these two flow components. Linear theory 
(see the discussion in \cite{Cb}) provides the dispersion relation
\begin{equation}\label{dis}
c=u_0-\frac{\gamma \tanh(kd)}{2k} \pm \frac{1}{2k}\,\sqrt{\gamma^2 \tanh^2(kd)+ 4gk\,\tanh(kd)}
\end{equation}
for the wave speed $c$, expressed in terms of the surface current speed $u_0$ in a flow of constant 
vorticity $\gamma$ in water of mean depth $d$; here $k=\dfrac{2\pi}{L}$ is the frequency of 
a wave with wavelength $L$. Note that (\ref{dis}) permits us to understand the effect of vorticity 
on the propagation speed of waves of small amplitude (that are realistically described by linear 
theory). For example, from the point of view 
of a fixed observer noticing right-propagating waves that interact with a current with vanishing 
surface speed (a setting corresponding to $c>0$ and $u_0=0$), from (\ref{dis}) we 
infer that the propagation speed is 
\begin{equation}\label{disv}
c(\gamma)=-\frac{\gamma \tanh(kd)}{2k} + \frac{1}{2k}\,\sqrt{\gamma^2 \tanh^2(kd)+ 4gk\,\tanh(kd)}\,,
\end{equation}
irrespective of the sign of $\gamma$. In this setting the underlying current beneath the flat 
free surface $y=0$ and above the bed $y=-d$ is given by $u(y)=\gamma y$, so that $\gamma>0$ 
corresponds to an adverse current since everywhere beneath the surface the current velocity opposes 
the direction of wave propagation, $\gamma<0$ corresponds to a favourable current (the current velocity 
points everywhere beneath the surface in the direction of wave propagation), while the 
irrotational case $\gamma=0$ is characterized by the absence of an underlying current. A 
simple analysis of (\ref{disv}) confirms that 
$$c(-|\gamma|)>c(0)>c(|\gamma|)$$
for every $\gamma \neq 0$. This means that favourable currents 
enhance the wave speed, while adverse currents reduce it. In contrast to 
this analysis of the wave speed, the first-order 
approximation that is pursued within the framework of linear theory is not conclusive in regard to 
the effect of underlying currents on the wave profiles: at first order, periodic travelling waves in 
flows of constant vorticity are sinusoidal (see \cite{Cb}). Consequently, nonlinear effects have 
to be accounted for and higher-order approximations are needed. We present in this paper an 
approach that provides the second and third order asymptotic expansion of waves of small 
amplitude, expressed in terms of a suitable amplitude parameter $b$. Following this procedure, 
we can readily obtain other characteristics of the water flow, such as the velocity field and the 
pressure beneath the wave. The knowledge of accurate approximations for these 
wave characteristics allows us to compare the results with those for irrotational flows. 
In the last part of the work we illustrate these flow characteristics for 
several different types of wave-current interactions. Laboratory experiments and numerical 
simulations for irrotational waves are discussed in \cite{Chen, Cl1, U}, while this type of studies for wave-current 
interactions in flows of constant vorticity were pursued in \cite{CKS1, DP, KS1, KS2, SCJ, T}. The results 
presented below permit a more detailed analysis and numerical simulation.

\section{Mathematical Formulation}

In this paper we discuss two dimensional, periodic travelling waves with constant vorticity. 
We assume that the water is incompressible and inviscid, over a flat bed and acted upon by gravity $g$. In what follows we make no shallowness or small amplitude approximation. In order to avoid extended revision we refer to the recent works \cite{Cb, CKS1, CS} for the mathematical formulation of the problem. In this section we make only a brief review of the problem, in order to provide sufficient terminology for the fundamental Boundary Value Problem (BVP) \eqref{op_basic}.

We are studying two-dimensional waves travelling at constant speed $c$. This means that 
in a two-dimensional frame moving with the constant speed $c$, the flow pattern --- and, 
in particular, the shape of the 
surface of the fluid --- does not change over time.

These assumptions, see \cite{CKS1}, allow the definitions of the free surface 
profile by
\begin{equation*}S = \set{(x,y): -\pi < x < \pi \text{ and }  y = 
\eta(x)}\,,\end{equation*} 
and the flat bed by
\begin{equation*}B = \set{(x,y): -\pi < x < \pi \text{ and } 
y=-d}\,,\end{equation*}
with $d>0$, for the normalized wavelength $2\pi$.

Without loss of generality, the assumption that the waves oscillate around the flat free surface $y=0$ is made, i.e., 
\begin{equation} \label{etn}
\int_{-\pi}^{\pi} \eta(x) dx = 0 .
\end{equation}

 We denote the velocity field of the flow by $(u(x,y),v(x,y))$, with 
$(x,y)\in\mathcal{D},$ where
\begin{equation*}
\mathcal{D} = \set{(x,y): -\pi < x < \pi \text{ and } -d < y < 
\eta(x)}\,.\end{equation*} 

In this paper we follow the assumption made in \cite{CS}, that there are no 
flow-reversals, which is formulated as the condition
\begin{equation}\label{ns}
u<c \quad\hbox{throughout the fluid}.
\end{equation}

We denote the pressure in the fluid by $P(x,y)$, with $(x,y)\in\mathcal{D}$. 
Neglecting the effects of surface tension --- a hypothesis that is appropriate for 
waves of moderate and large amplitude (see the discussion in \cite{Cb}) --- we 
impose that the water pressure is constant on $S$, 
$P(x,\eta(x))=P_{atm}$, where $P_{atm}$ is the atmospheric pressure.

Let us recall the following definitions:
\begin{itemize}
\item The vorticity of the flow is defined by 
\begin{equation}\label{gamma}
\gamma := u_y-v_x ,
\end{equation} 
which for 
the rest of the paper we will assume constant. 

\item The relative mass flux\footnote{The terminology `relative' is due to the 
fact that $(u-c)$ is the relative horizontal velocity of the flow, with 
reference to the moving frame at speed $c$.} is defined by
\begin{equation}\label{mf}
p_0:=\int_{-d}^{\eta(x)}\big(u(x,y)-c\big)\,dy<0,
\end{equation} 
which in fact is independent of $x$, see \cite{CKS1}. Moreover, \eqref{mf} and 
\eqref{ns} show that $p_0$ is negative. This relation shows that the amount of 
water passing any vertical line is constant throughout $\mathcal{D}$.

\item The stream function $\psi(x,y)$ is defined as the unique solution of the 
differential equations  
\begin{equation}
 \label{ds}
 \psi_x = -v, \qquad \psi_y = u-c \text{ in }\  \overline{\mathcal{D}}\,,
\end{equation}
subject to 
\begin{equation}
 \label{bs}
 \psi(x,-d) = -p_0\;.
\end{equation}

\item In \cite{CS} it was proven that 
the expression 
\begin{equation}
\label{eq:E}
 \frac{(u-c)^2+v^2}{2}+ gy+ P -\gamma\psi
\end{equation}
equals a constant $E$ throughout $\mathcal{D}$. 
The constant $Q=E-P_{atm}$ is called the hydraulic head.

\end{itemize}

We also recall some basic terminology, as used in \cite{CS}:
\begin{itemize}
\item Firstly, the Dubreil-Jacotin transformation maps the unknown domain 
$\mathcal{D}$ to the rectangle 
\begin{equation}
 R=\set{(q,p):\ -\pi<q<\pi\,,\  p_0<p<0}\;,
\end{equation}
\begin{figure}[ht!]
\centering
\includegraphics[height = 45mm, width =95mm]{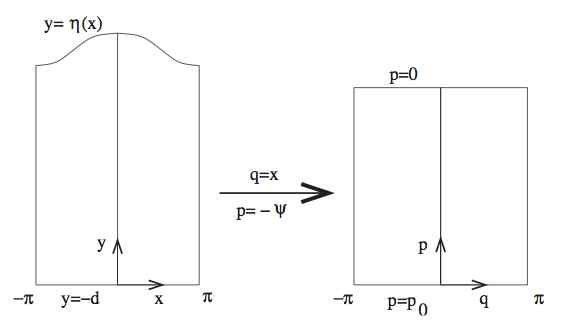}
\caption{Dubreil-Jacotin transformation}
\label{D-J}
\end{figure} 
The independent variables as they appear in the Dubreil-Jacotin transformation 
\cite{DJ} are
\begin{equation*}
q=x, \qquad p=-\psi.
\end{equation*}
\item Secondly, the height function of the wave above $B$ is defined as 
\begin{equation} \label{def:h} h(q,p) = y + d\;. \end{equation} 
This function satisfies the following nonlinear BVP on the rectangle $R$. 
We omit here the details of this construction, and refer to \cite{CS}. 
\end{itemize}

\begin{definition}
The constitutive equations for the height function $h(q,p)$, which 
is even and $2\pi$-periodic in $q$, are the following:
\begin{equation}\label{op_basic}
\left\{\begin{array}{l}
\mathcal{H}[h] := (1+h_q^2)h_{pp} - 2h_ph_qh_{pq} +h_p^2h_{qq} -\gamma h_p^3 
=0\ 
 \text{ on } \ R\,,\\
\mathcal{B}_0[h] := 1+h_q^2(q,0) + (2gh-Q)h_p^2(q,0) = 0 \,, \\
\mathcal{B}_1[h] := h(q,p_0) = 0\,.
\end{array}\right.
\end{equation}

Moreover, the free boundary $S$ is given by the 
expression $$h(q,0)=\eta(x)+d.$$
\end{definition}
In this equation $\gamma$ denotes the vorticity as defined in \eqref{gamma}, which 
is assumed constant, being indicative of underlying 
currents. Vanishing vorticity is the property typical of uniform currents and a 
constant 
vorticity characterizes the linearly sheared tidal currents. 
We also recall that $g$ is the gravitational constant  of acceleration, $Q$ is the hydraulic head 
and $p_0$ is the relative mass flux, introduced above.
Moreover, we define the height of the wave as the maximal variation of the 
oscillations of the free surface, given by 
\begin{equation}
 \label{eq:amplitude}
a:=\max_{q\in [-\pi,\pi]} h(q,0) - \min_{q\in [-\pi,\pi]} h(q,0) = 
h(0,0)-h(\pi,0).
\end{equation}

In fact, the BVP \eqref{op_basic} is a reformulation of the Euler equations 
restricted to $\mathcal{D}$, under the assumptions made previously in this 
section.

In the following subsections we discuss some results from \cite{CS}, and provide 
a new point of view, which allows a generalization of these results.

\subsection{Laminar solutions}

The laminar flows are readily obtained as the $q$-independent solutions of the boundary-value 
problem (\ref{op_basic}), given by the following formula
\begin{equation} \label{H}
H(p;\lambda) =  \frac{2(p-p_0)}{\sqrt{\lambda - 2 \gamma p}+\sqrt{\lambda-2\gamma p_0}}\,,\quad p_0 \le p \le 0\,,
\end{equation}
provided that the parameter $\lambda>0$ satisfies the equation
\begin{equation} \label{Q}
Q =\lambda + \frac{4 g |p_0|}{\sqrt{\lambda} + \sqrt{\lambda- 2 \gamma p_0}}\,.
\end{equation} 

\subsection{Solutions of the linearised problem}

Apart from the laminar flows one can obtain the solution of the linearised problem. For this, 
the BVP \eqref{op_basic} is linearised around the laminar flow $H(p;\lambda)$. For the  specific value $\lambda_*$ that satisfies the dispersion equation 
\begin{equation}\label{disp}
\dfrac{ \lambda }{g- \gamma \sqrt{ \lambda }} + \tanh \left( \dfrac{2p_0}{ \sqrt{ \lambda } +\sqrt{ \lambda -2 p_0 \gamma } }\right) =0 ,
\end{equation}
 an existence result for the solution of this BVP is presented in \cite{CS}. Moreover, the linearised solution is given by 
$\hat{h}(q,p;b)= H(p;\lambda_*) + b \ m(q,p)$, where
\begin{equation}\label{sol-lin}
m(q,p) =\dfrac{\sqrt{ \lambda_* -2 p_0 \gamma }}{\sqrt{ \lambda_* -2 p \gamma }}\sinh\left(  \frac{2(p-p_0)}{\sqrt{\lambda_* - 2 \gamma p}+\sqrt{\lambda_*-2\gamma p_0}}\right) \  \cos q .
\end{equation}
In our investigation we start from this result and we obtain an extension of it by making a two-fold interpretation.

Firstly, we interpret the function 
\begin{equation}
\label{sol_perturb}
\hat{h}(q,p;b)= H(p;\lambda_*) + b \ m(q,p), 
\end{equation}
as a perturbation of a laminar solution of the system  \eqref{op_basic}, in the following sense: the system \eqref{op_basic} is satisfied up to order $b^2$, i.e.,
\begin{equation}
\label{eq:ob2}
 \mathcal{H}[\hat{h}](q,p) = \mathcal{O}(b^2)\,,\; \mathcal{B}_0[\hat{h}](q) = \mathcal{O}(b^2) \text{ and } \mathcal{B}_1[\hat{h}](q) = 0\;
 \end{equation}
 and the height of the water wave is of order $b$, i.e, 
$$\hat{h}(0,0;b)-\hat{h}(\pi,0;b)=b\left[ m(0,0)-m(\pi,0)\right]=\mathcal{O}(b). $$
 Having in mind that the wave height vanishes for laminar flows, we can view the expression \eqref{sol_perturb} as an approximation of small amplitude water waves. 
 
Secondly, we can regard this expression as the first order asymptotic expansion of the exact solution to the BVP \eqref{op_basic}, 
which leads to the natural question of finding higher-order terms. This point of view is 
motivated by a bifurcation argument, provided in \cite{CS} and discussed below. 
 
\subsection{Bifurcation}
 
We define the curve that represents the laminar flows $$\mathcal{T}=\left\{ ( 
Q(\lambda), H(p;\lambda)) : \lambda>0 \right\},$$   with $Q$ and $H$ given by 
\eqref{Q} and \eqref{H}, respectively.
 
\begin{figure}[ht!]
  \centering
  \includegraphics[height = 40mm, width =80mm]{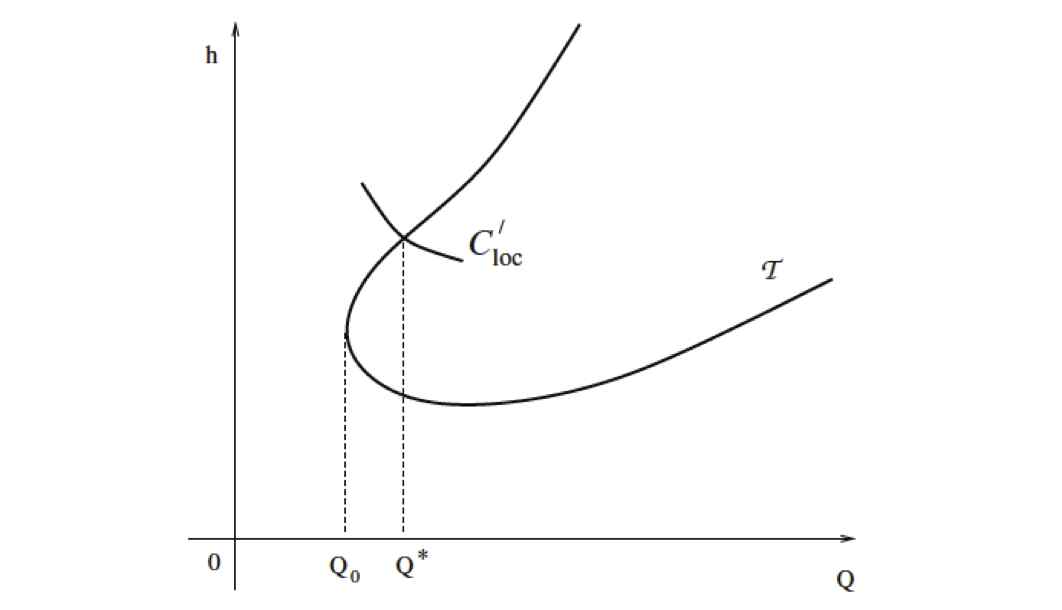}
  \caption{Bifurcation on the curve of laminar flows.}
  \label{fig:bifurc}
\end{figure}

Restricting our attention to rotational flows with constant vorticity $\gamma$, 
in \cite{CS} it was proven that near this curve, as the parameter $\lambda$ 
varies, 
there are generally no genuine waves, except at critical values 
$\lambda=\lambda_*$ determined by the dispersion relation \eqref{disp}. Near 
this bifurcating laminar flow $H^\ast$, we have two solution curves: one 
laminar 
solution curve $\lambda \mapsto H(p;\lambda)$, where $\lambda$ and $Q$ are 
related by \eqref{Q}, and one non-laminar solution curve $Q \mapsto h(q,p;Q)$ 
such that $h_q \not \equiv 0$ unless $h=H^\ast$, see Figure \ref{fig:bifurc}. 

\subsection{Higher order approximation}

Using the above argument we observe that the variation of the parameter $Q$ 
(which is the hydraulic head of the flow) about the uniquely determined 
value $Q^*$ implies an approximation of non-laminar flows, in the following 
sense. Consider the approximation for the hydraulic head of the flow,
$$Q 
\approx Q^{(m)}=Q(b)=Q^*+\sum_{n=1}^m Q_n b^n,$$
and the approximation for the height function $h(q,p;Q)$, 
\begin{equation}\label{h_appr_aux}h(q,p;Q)\approx h^{(m)}(q,p;b)=\sum_{n=0}^m 
h_n(q,p;Q^*) b^n,\end{equation} 
with $h_0(q,p;Q^*)\equiv H(p,\lambda_*)$.

Our goal is to determine $\{ h_n(q,p;Q^*)\}_{n=1}^m$ such that the 
system \eqref{op_basic} is satisfied up to order $b^{m+1}$, i.e.,
\begin{equation}
\label{eq:obm}
 \mathcal{H}[h^{(m)}](q,p) = \mathcal{O}(b^{m+1})\,,\; \mathcal{B}_0[h^{(m)}](q) 
= 
\mathcal{O}(b^{m+1}) \text{ and } \mathcal{B}_1[h^{(m)}](q) = 0\; .
 \end{equation}

Moreover, the dominant term of the wave height of this approximation of the water 
wave is of order $b$, i.e, 
\begin{equation}
\label{amp:obm}a^{(m)}=h^{(m)}(0,0;b)-h^{(m)}(\pi,0;b)=b \left[ h_1(0,0;Q^*) - 
h_1(\pi,0;Q^*) \right] + \mathcal{O}(b^2).
\end{equation}
This condition indicates the accuracy of approximation of genuine waves in the following 
sense: The relation \eqref{eq:obm} shows that $h^{(m)}$ is an approximation of a 
solution up to order $b^m$, 
whereas the relation \eqref{amp:obm} shows that $h^{(m)}$ differs from the 
laminar solution at order $b$. Thus $h^{(m)}$ is `closer' to a non-laminar 
solution.

In Section 3, we present the results for $m=2$. In Section 4, where we work out  
the case $m=3$, the situation is different. There we observe that there is no $ 
h^{(3)}(q,p;b)$ such that {\it {both}}  conditions
 $\displaystyle \mathcal{H}[h^{(3)}](q,p) = \mathcal{O}(b^4)\, \text{ and }  
\mathcal{B}_0[h^{(3)}](q) = \mathcal{O}(b^{4})$ are satisfied. This 
means that the approximation made on \eqref{h_appr_aux} is accurate up to $m=2$, 
and there appears an additional error for $m\geq 3$. Thus, we introduce an auxiliary 
parameter which allows the control of this additional error. Indeed, we 
construct the following estimates
$$ \left\| \mathcal{H}[h^{(3)}](q,p) \right\| \leq \epsilon_1 b^3\, \text{ and } 
\left\| \mathcal{B}_0[h^{(3)}](q)\right\| \leq \epsilon_2 b^3,$$
for some small $\epsilon_1\geq 0$ and $\epsilon_2\geq 0$, which depend on this 
auxiliary parameter.  

Moreover, for particular values of the auxiliary parameter we can make each of 
$\epsilon_1$ and $\epsilon_2$ arbitrarily small (but not simultaneously), i.e.,
\begin{itemize}
\item  either, $\displaystyle \mathcal{H}[h^{(3)}](q,p) = \mathcal{O}(b^4)$,  
\item or,  $\displaystyle\mathcal{B}_0[h^{(3)}](q) = \mathcal{O}(b^{4})$.
\end{itemize}
The selection of the proper value of the auxiliary parameter is discussed in the 
last parts of Section 4.1 and Section 4.2, as well as in Section 5.2.

\section{Second Order Asymptotics}

The problem of the extension of the asymptotic expansion of the solution of \eqref{op_basic} takes the following form: 
with $\lambda_*$ defined by \eqref{disp}$, H(p;\lambda_*)$ given by \eqref{H} and $m(q,p)$ given by \eqref{sol-lin}, 
determine $u(q,p)$ such that the function 
\begin{equation}
\label{sol_perturb_2}
h^*(q,p;b)= H(p;\lambda_*) + b \ m(q,p) + b^2 u(q,p), 
\end{equation}
satisfies the system \eqref{op_basic} up to the second order, i.e.
\begin{equation}
\label{eq:ob3}
 \mathcal{H}[h^*](q,p) = \mathcal{O}(b^3)\,,\; \mathcal{B}_0[h^*](q) = \mathcal{O}(b^3) \text{ and } \mathcal{B}_1[h^*](q) = 0.\;
 \end{equation}
Note that, in the absence of stagnation points, the real-analyticity of the wave profile for flows of constant vorticity, 
established in \cite{CE}, ensures the existence of a power series expansion. Since the first-order (linear) approximation 
of the wave profile fails to capture the effects of the underlying current, the issue of obtaining high-order 
approximations should be explored in this regard.

\subsection{Irrotational flow}

For the purpose of comparison we analyse first the irrotational case $\gamma=0$. 
This simpler case provides us with the opportunity to explain in more detail 
the strategy we use throughout this paper to obtain higher order terms of the expansion of the solution.

In the irrotational case, the laminar solution $H(p;\lambda_*)$ simplifies to 
\begin{equation}\label{H*0}
H(p;\lambda_*)=\dfrac{p-p_0}{\sqrt{ \lambda_*}}\,,
\end{equation}
the linearised solution $m(q,p)$ becomes
\begin{equation}\label{m0}
m(q,p) = \sinh\left( \dfrac{p-p_0}{\sqrt{ \lambda_* }}\right) \ \cos q,
\end{equation}
where $\lambda_*>0$ satisfies, in this setting, the dispersion relation
\begin{equation}\label{disp0}
\lambda + g \tanh \left( \dfrac{p_0}{\sqrt{ \lambda }}\right) = 0\,.
\end{equation}
Moreover the corresponding value of $Q$ is 
\begin{equation}\label{Q*}Q^* = \lambda_* - \dfrac{2 g p_0}{\sqrt{ \lambda_*}}\;.\end{equation}

We first apply the operator $\mathcal{H}$, defined in \eqref{eq:ob3}, to the expression \eqref{sol_perturb_2}. The outcome will be a polynomial in the $b$ variable, of order 
up to $b^6$. The specific above choices yield the vanishing of the zero and first order of the polynomial; for a detailed discussion of this aspect, we refer to the \cite{CKS1}. Our present goal 
is the vanishing of the $b^2$ coefficient of the polynomial. 

In particular, using the fact that $H_q = 0$, we get 
\begin{equation*}
\begin{split}
\mathcal{H}[h^*] =  H_{pp} &+ b \big( m_{pp} + H_p^2 m_{qq} \big) \\
& + b^2 \left[ u_{pp}+ H_p^2 u_{qq} + H_{pp} - 2H_p (m_qm_{qp}-m_pm_{qq})    \right]  + \mathcal{O}(b^3).
\end{split}
\end{equation*}
Observing that \eqref{H*0} yields \begin{equation} \label{H0_help}
H_p \equiv \frac{1}{\sqrt{\lambda_*}} \ \text{ and } \ H_{pp} \equiv 0  \qquad \text{for all}\quad p\in [p_0,0]
\end{equation} and substituting\eqref{m0} in the previous expression for $\mathcal{H}[h^*]$, we get
\begin{equation}
\mathcal{H}[h^*] = b^2 \left[ u_{pp}+\dfrac{1}{\lambda_*} u_{qq} -\dfrac{1}{\lambda_*}\sinh\left(2\frac{p-p_0}{\sqrt{\lambda_*}}\right)   \right]  + \mathcal{O}(b^3).
\end{equation}
Proceeding similarly for the boundary condition, i.e. for the operator $\mathcal{B}_0[h^*]$ in \eqref{eq:ob3}, yields 
\begin{equation*}
 \begin{aligned}
\mathcal{B}_0[h^*] &=1 + (2gH-Q)H_p^2  + b \left\{ 2H_p \big[(2gH-Q) m_{p}+ g H_p m \big] \right\}  \\ 
&+ b^2 \left[ m_q^2 + 4 g m H_p m_p + (2gH-Q)(m_p^2 + 2 H_p u_p) +2g H_p^2 u \right]\\ &+ \mathcal{O}(b^3)\;,
 \end{aligned}
\end{equation*}
evaluated at $p=0$. Applying the equalities \eqref{H*0}, \eqref{m0}, \eqref{Q*} and \eqref{H0_help} in the previous expression and using the fact that \eqref{disp0} is equivalent to the identities 
\begin{equation}\label{disp_simpl}
\cosh\left( \dfrac{p_0}{\sqrt{ \lambda_* }}\right) = \dfrac{g}{\sqrt{g^2-\lambda_*^2}} \quad \text{ and } \quad \sinh\left( \dfrac{p_0}{\sqrt{ \lambda_* }}\right) = -\dfrac{\lambda_*}{\sqrt{g^2-\lambda_*^2}}
\end{equation} 
we get that
\begin{equation}
\begin{split}
\mathcal{B}_0[h^*] = -b^2 2\sqrt{\lambda_*}  & \left\{  \left(u_{p} + \frac{g}{\lambda_*^{3/2}}u\right)\bigg|_{p=0} \right. \\ & \left. - \frac{1}{4\sqrt{\lambda_*}}\dfrac{1}{g^2-\lambda_*^2}\bigg[ 3g^2+\lambda_*^2 + \big( 3g^2-\lambda_*^2\big) \cos(2q)\bigg]\right\}  \\ & + \mathcal{O}(b^3).
\end{split}
\end{equation}

Consequently, we obtain the constitutive equations for the function $u(q,p)$, even and $2\pi$ periodic in $q$ variable,
\begin{equation}\label{eq:u}
\begin{aligned}
u_{pp}+\dfrac{1}{\lambda_*} u_{qq} &= \dfrac{1}{\lambda_*}\sinh\left(2\frac{p-p_0}{\sqrt{\lambda_*}} \right) , \qquad &(q,p) \in R, \\
 u_{p} - \frac{g}{\lambda_*^{3/2}}u  &= \frac{1}{4\sqrt{\lambda_*}}\dfrac{1}{g^2-\lambda_*^2}\bigg[ 3g^2+\lambda_*^2 + \big( 3g^2-\lambda_*^2\big) \cos(2q)\bigg], &p=0, \\
u&=0, \qquad &p=p_0.
\end{aligned}
\end{equation}
In order to solve the above BVP we use separation of variables and we take in consideration the following facts:
\begin{itemize}
\item A special solution of the PDE is given by $ \displaystyle\frac{1}{4}\sinh\left(2\frac{p-p_0}{\sqrt{\lambda_*}} \right).$
\item The function $u(q,p)$ is even and $2\pi$ periodic.
\item We must have $u(q,p_0)=0$.
\end{itemize}
Therefore, we get the general form
\begin{equation}\label{u}u(q,p)=\frac{1}{4}\sinh\left(\frac{2(p-p_0)}{\sqrt{\lambda_*}} \right) +
 A_0 (p-p_0) + \sum_{n=2}^\infty A_n \cos(nq) \sinh\left(\frac{n(p-p_0)}{\sqrt{\lambda_*}} \right)
\end{equation}
and we have to determine the constants $A_0$ and $\{A_n\}_2^\infty$ such that the Robin Boundary Condition (BC) at $p=0$ is satisfied. Indeed, the separation of variables says that the general solution of the 
Partial Differential Equation 
$$ u_{pp}+\dfrac{1}{\lambda_*} u_{qq}=0$$ on the rectangle is a bilinear 
combination of trigonometric functions in $q$ and exponential (hyperbolic) 
functions in $p$ i.e. $\{\cos(nq),\sin(nq)\}, \ n\in \N$ and $\exp(\pm n 
p/\sqrt{\lambda_*}), \ n\in \N$. This is due to the periodicity in 
$q$.\footnote{Once we have obtained $\{\cos(nq),\sin(nq)\}$, separation 
of variables gives the exponentials $\exp(\pm n p/\sqrt{\lambda_*})$.} The 
evenness in $q$ excludes the dependence on $\sin(nq), \ n\in \N$ and the 
condition $u=0$, on $p=p_0$ gives exactly the 
$\sinh\left(n\frac{p-p_0}{\sqrt{\lambda_*}} \right)$ dependence.

Thus, applying the
form \eqref{u} to the second of the equations \eqref{eq:u}, we get
$$ A_0=\frac{\lambda_*}{4} \dfrac{3g^2-\lambda_*^2}{g^2-\lambda_*^2} \frac{1}{g 
p_0+ \lambda_*^{3/2}}, \qquad A_2=\frac{3g^2-\lambda_*^2}{8\lambda_*^2} \quad 
\text{ and } \quad A_n=0, \  n\geq3,$$
with $\lambda_*$ given as the solution of \eqref{disp0}.

\subsection{Constant Vorticity}

The strategy 
is conceptually the same for the rotational case when the vorticity $\gamma$ is 
constant. The first goal is to derive the analogue of the BVP \eqref{eq:u}. In 
this case, the laminar solution $H(p)$ is given by
\begin{equation}\label{H*0_r}
H(p)=2\dfrac{p-p_0}{ r(p)+r(p_0)}\,,
\end{equation}
where \begin{equation}\label{r(p)}
r(p)=\sqrt{\lambda_*-2\gamma p},
\end{equation}
with $\lambda_*>0$ satisfying the dispersion relation \eqref{disp}. Moreover, 
the linearised solution $m(q,p)$ is given by \eqref{sol-lin} and rewritten as 
follows 
\begin{equation}\label{m0_r}
m(q,p) = \dfrac{r(p_0)}{r(p)}\sinh\big( H(p)\big) \ \cos q.
\end{equation}

The corresponding value of $Q$, using \eqref{Q}, is 
\begin{equation}\label{Q*_r}Q^* = \lambda_* - \dfrac{4 g p_0}{\sqrt{ \lambda_*}+\sqrt{\lambda_*-2\gamma p_0}}\;.\end{equation}

Then we evaluate the system \eqref{op_basic} for $h^*$, defined by
\begin{equation}
\label{sol_perturb_2_r}
h^*(q,p;b)= H(p) + b \ m(q,p) + b^2 u(q,p), 
\end{equation}
with $H(p)$ and $m(q,p)$ given by \eqref{H*0_r} and \eqref{m0_r}, respectively. 
Demanding the validation of this system 
up to the second order, i.e.
\begin{equation}
\label{eq:ob3_r}
 \mathcal{H}[h^*](q,p) = \mathcal{O}(b^3)\,,\; \mathcal{B}_0[h^*](q) = 
\mathcal{O}(b^3) \text{ and } \mathcal{B}_1[h^*](q) = 0,\;
 \end{equation}
yields the analogue of the BVP \eqref{eq:u}. After tedious but straightforward 
calculations we get the following form for the BVP, for the even and 
$2\pi$-periodic function $u(q,p)$ on the rectangle  
$R\equiv[-\pi,\pi]\times[p_0,0]$:
\begin{equation}\label{eq:u_r}
\begin{aligned}
u_{pp}(q,p)+\dfrac{1}{r^2(p)} u_{qq}(q,p) -\frac{3\gamma}{r^2(p)}u_p(q,p)&= 
f_0(p)+f_2(p) \cos(2q) , \\
 u_{p}(q,0) - \frac{g}{\lambda_*^{3/2}}u(q,0)  &=  g_0+g_2 \cos(2q), \\
u(q,p_0)&=0,
\end{aligned}
\end{equation}
where
\begin{equation*}\label{f0}
\begin{aligned}
f_0(p) = -\dfrac{r^2(p_0)}{r^7(p)}  \bigg[ 3\gamma^3  &- 3\gamma \left( 
\gamma^2+2r^2(p) \right) \cosh\left( 2 \frac{r(p)-r(p_0)}{\gamma} \right) \\ & + 
2 r(p) 
\left( 3\gamma^2+2r^2(p) \right) \sinh\left( 2 \frac{r(p)-r(p_0)}{\gamma}\right) 
 \bigg],
\end{aligned}
\end{equation*}
\begin{equation*}\label{f2}
\begin{aligned}
f_2(p) = -\frac{\gamma}{4}\dfrac{r^2(p_0)}{r^7(p)}  \bigg[ 3\gamma^2 -2 r^2(p)  
&-  \left( 
3\gamma^2+4r^2(p) \right) \cosh\left( 2 \frac{r(p)-r(p_0)}{\gamma} \right) \\ & 
+ 6 
\gamma r(p) \sinh\left( 2 \frac{r(p)-r(p_0)}{\gamma}\right)  \bigg],
\end{aligned}
\end{equation*}
\begin{equation*}\label{g0}
\begin{aligned}
g_0 = \dfrac{r^2(p_0)}{2\lambda_*}  
\dfrac{3g^2+\lambda_*^2}{g^2-\lambda_*^2-\gamma\sqrt{\lambda_*}\left( 
2g-\gamma\sqrt{\lambda_*}\right) },
\end{aligned}
\end{equation*}
\begin{equation*}\label{g2}
\begin{aligned}
g_2 = \dfrac{r^2(p_0)}{2\lambda_*}  
\dfrac{3g^2-\lambda_*^2}{g^2-\lambda_*^2-\gamma\sqrt{\lambda_*}\left( 
2g-\gamma\sqrt{\lambda_*}\right) }.
\end{aligned}
\end{equation*}

\begin{remark}
The above BVP \eqref{eq:u_r} reduces to \eqref{eq:u} for $\gamma=0$, i.e. 
the function 
$f_2(p)$ vanishes when $\gamma$ vanishes.
\end{remark}

Applying separation of variables to \eqref{eq:u_r}, we obtain the following 
formulation of the above BVP: Determine $u(q,p)$ given by
\begin{equation} \label{u_b2}
u(q,p)=u_0(p)+ \sum_{n=2}^\infty u_n(p) \cos(nq),
\end{equation}
where the functions $u_n(p)$ satisfy the following BVPs:
\begin{equation}\label{eq:u_r_n}
\begin{aligned}
u''_n(p)-\dfrac{n^2}{r^2(p)} u_n(p) -\frac{3\gamma}{r^2(p)}u'_n(p)&= f_n(p) , 
\qquad p\in(p_0,0)\\
 u'_n(0) - \frac{g}{\lambda_*^{3/2}}u_n(0)  &=  g_n, \\
u_n(p_0)&=0,
\end{aligned}
\end{equation}
where $f_0(p)$, $f_2(p)$, $g_0$ and $g_2$ are defined above. Moreover, 
$f_n\equiv0$ and $g_n=0$, for $n\geq3$.

Now, we make the following change of variables
\begin{equation} \label{transform}
u_n(p)=\dfrac{\hat{u}_n(r(p))}{r(p)}, \qquad \text{ with } \ 
r(p)=\sqrt{\lambda_*+2\gamma p}
\end{equation}
which transforms the above BVPs into the following:
\begin{equation}\label{eq:u_r_n-hat}
\begin{aligned}
\dfrac{d^2}{dr^2}\hat{u}_n(r)-\dfrac{n^2}{\gamma^2} \hat{u}_n(r) &= 
\hat{f}_n(r) 
, \qquad &r\in\left(r_0,\sqrt{\lambda_*}\right)\\
 \dfrac{d}{dr}\hat{u}_n - \dfrac{g + \gamma 
\sqrt{\lambda_*}}{\gamma\lambda_*}\hat{u}_n  &=  \hat{g}_n, \qquad 
&r=\sqrt{\lambda_*}\\
\hat{u}_n&=0, \qquad &r=r_0:=r(p_0)=\sqrt{\lambda_*+2\gamma p_0}.
\end{aligned}
\end{equation}

Here $\hat{f}_n$ and $\hat{g}_n$ are defined by
$$\hat{f}_n(r) = \hat{f}_n(r(p)) = \frac{r^3(p)}{\gamma^2} f_n(p) \qquad \text{ 
and } \qquad \hat{g}_n=-\frac{\lambda_*}{\gamma}g_n, \qquad n\in \N.$$

The solution basis of the homogeneous differential equations appearing in BVP 
\eqref{eq:u_r_n-hat} is given by 
\begin{itemize}
\item $\{1,r\}$ for $n=0$ and
\item $\left\{ \sinh\left( n \dfrac{r-r_0}{\gamma} \right),\cosh\left( n 
\dfrac{r-r_0}{\gamma} \right) \right\}$ for $ n\geq 1.$
\end{itemize}

Using the method of variation of parameters, we obtain the special solutions of 
the BVP \eqref{eq:u_r_n-hat}. In particular, we seek for special solutions of 
the inhomogeneous problems of the form
\begin{itemize}
\item For $n=0$: $$\hat{u}_0^{(s)}(r)=\alpha_0(r) + \beta_0(r) r,$$ under the 
constraint  $$\alpha'_0(r) + \beta'_0(r) r = 0.$$ 
\item For $n\geq 1$: $$\hat{u}_n^{(s)}(r)=\alpha_n(r) \sinh\left( n 
\dfrac{r-r_0}{\gamma} \right) + \beta_n(r) \cosh\left( n \dfrac{r-r_0}{\gamma} 
\right),$$  under the constraint  $$\alpha'_n(r) \sinh\left( n 
\dfrac{r-r_0}{\gamma} \right) + \beta'_n(r) \cosh\left( n \dfrac{r-r_0}{\gamma} 
\right)=0.$$
\end{itemize}
The insertion of these forms in the inhomogeneous differential equations 
appearing in BVP \eqref{eq:u_r_n-hat} yields $\{\alpha_n,\beta_n\}, \ n \in \N$. 
Consequently, the special solutions read 

$$\hat{u}_0^{(s)}(r)=-\gamma r_0^2\dfrac{1-\cosh\left( 2 \dfrac{r-r_0}{\gamma} 
\right)}{8r^2}- r_0^2\dfrac{\sinh\left( 2 \dfrac{r-r_0}{\gamma} \right)}{4r}$$ 
and
$$\hat{u}_2^{(s)}(r)=-\gamma r_0^2\dfrac{1-\cosh\left( 2 \dfrac{r-r_0}{\gamma} 
\right)}{8r^2}- r_0^2\dfrac{\sinh\left( 2 \dfrac{r-r_0}{\gamma} \right)}{4r}.$$

The next step is to apply the linear combination $$\hat{u}_n = \hat{u}_n^{(s)} 
+ 
E_n\hat{u}_n^{(g)}$$  to the Robin condition of \eqref{eq:u_r_n-hat} and obtain 
the specific value for $E_n$. One can readily see that $E_n=0$ for $n\geq 3$.

Finally, after determining $E_0$ and $E_2$, we use \eqref{transform} and, by  
means of \eqref{u_b2}, we obtain the following formula for the second order term of the 
expansion \eqref{sol_perturb_2_r}:
\begin{equation}\label{u_r}
\begin{aligned}
u(q,p)&= r^2(p_0) \left[ C_0 \dfrac{H(p)}{r(p)} -\gamma \dfrac{1-\cosh\left( 2 H(p) \right)}{8r^3(p)}+ \dfrac{\sinh\left( 2 H(p) \right)}{4r^2(p)} \right] \\ &+r^2(p_0) \cos(2q) \left[  C_2 \dfrac{\sinh\left( 2 H(p) \right)}{r(p)} -\gamma \dfrac{1-\cosh\left( 2 H(p) \right)}{8r^3(p)}+ \dfrac{\sinh\left( 2 H(p) \right)}{4r^2(p)} \right],
\end{aligned}
\end{equation}
where 
$$C_0=\frac{\sqrt{\lambda_*}\left[  3g \left(g -\gamma\sqrt{\lambda_*} \right) + \lambda_* \left( \gamma^2-\lambda_*\right) \right] \left( r(p_0) + \sqrt{\lambda_*}\right)  }{4\left[ \left(g -\gamma\sqrt{\lambda_*} \right)^2 -\lambda_*^2\right] \left[2gp_0 +\lambda_* r(p_0) + \sqrt{\lambda_*} r^2(p_0) \right] }$$
and 
$$C_2=\dfrac{3g\left(g -\gamma\sqrt{\lambda_*} \right) + \lambda_* \left( \gamma^2-3\lambda_*\right)  }{8\lambda_*^{5/2}}.$$

Therefore, the function $h^*$ given by \eqref{sol_perturb_2_r} satisfies the 
system \eqref{eq:ob3_r}. We note that for the value $\gamma=0$, the expression 
\eqref{u_r} takes the form \eqref{u}, meaning that we obtain the irrotational 
case as a special case of the constant vorticity formulation.

\section{Third Order}

With a similar methodology as that used to derive approximations of second order, 
one can derive approximation of third order.
The results are provided below; the reader can easily verify them 
with some Computer Algebra Software.

In particular, we use the formula 
\begin{equation}
\label{sol_perturb_3}
\tilde{h}(q,p;b)= H(p) + b \ m(q,p) + b^2 u(q,p) + b^3 w(q,p), 
\end{equation}
with $H$, $m$ and $u$ given in the previous section and $\lambda_*$ given by 
\eqref{disp}. Then our aim is to 
determine $w$ so that $\tilde{h}$ satisfies the system \eqref{op_basic} up to 
the third order, i.e.
\begin{equation}
\label{eq:ob4}
 \mathcal{H}[\tilde{h}](q,p) = \mathcal{O}(b^4)\,,\; 
\mathcal{B}_0[\tilde{h}](q) 
= \mathcal{O}(b^4) \text{ and } \mathcal{B}_1[\tilde{h}](q) = 0.\;
 \end{equation}
However, there is no such $\tilde{h}(q,p;b)$. Thus, instead, we introduce an 
auxiliary parameter $\tilde{B}$ (for the irrotational case we denote it by $B$) 
which allows the control of this additional error estimate. Indeed, we construct 
the following estimates
$$ \left\| \mathcal{H}[h^{(3)}](q,p) \right\| \leq \epsilon_1 b^3\, \text{ and } 
\left\| \mathcal{B}_0[h^{(3)}](q)\right\| \leq \epsilon_2 b^3,$$
for some small $\epsilon_1\geq 0$ and $\epsilon_2\geq 0$, which depend on this 
auxiliary parameter.

\subsection{Irrotational case}

Following the organisation of the previous section we first discuss the case 
that $\gamma=0$. Therefore, if we apply \eqref{sol_perturb_3} to the system 
\eqref{op_basic}, we get a polynomial on $b$ of finite order (in particular of 
order $b^9$). The requirement that the coefficient of the $b^3$ term vanishes 
yields the following BVP for the even and $2\pi$-periodic (in the $q$ variable) 
function $w(q,p)$:
\begin{equation}\label{eq:w}
\begin{aligned}
w_{pp}+\dfrac{1}{\lambda_*} w_{qq} &= f_1(p) \cos q + f_3(p)\cos(3q) , \qquad 
&(q,p) \in R, \\
 w_{p} - \frac{g}{\lambda_*^{3/2}}w  &= g_1 \cos q + g_3\cos(3q), &p=0, \\
w&=0, \qquad &p=p_0,
\end{aligned}
\end{equation}
where $f_1(p)$ and $f_3(p)$ are known functions (bilinear combinations of 
hyperbolic functions and polynomials) 
and $g_1$ and $g_3$ are known constants, dependent on the parameters $g, \ p_0$ 
and $\lambda_*$. We use the fact that this problem has the same general 
solution 
(for the homogeneous case) with the problem \eqref{eq:u} and we determine its 
special solutions by the method of variation of parameters.
Then the expansion \begin{equation*}
\tilde{h}(q,p;b)= H(p) + b \ m(q,p) + b^2 u(q,p) + b^3 \tilde{w}(q,p), 
\end{equation*}
with $\tilde{w}(q,p)$ given the following formula
\begin{equation}\label{w-ir}
\begin{aligned}
\tilde{w}(q,p)&=\left[  B_1 \left( p-p_0\right) \cosh \left( 
\dfrac{p-p_0}{\sqrt{\lambda_*}}\right) + B_2 \sinh \left( 
3\dfrac{p-p_0}{\sqrt{\lambda_*}}\right)\right] \cos q \\
&+\left[  D_1 \sinh \left( \dfrac{p-p_0}{\sqrt{\lambda_*}}\right) + D_2 \sinh 
\left( 3\dfrac{p-p_0}{\sqrt{\lambda_*}}\right)\right] \cos(3q) ,
\end{aligned}
\end{equation}
with $$B_1 = \frac{\lambda_*}{4} \dfrac{3g^2-\lambda_*^2}{g^2-\lambda_*^2} 
\frac{1}{g p_0+ \lambda_*^{3/2}}, \qquad B_2=\dfrac{9g^2}{32\lambda_*^2},$$
$$D_1 = -\dfrac{3g^2-\lambda_*^2}{32\lambda_*^2}  \qquad \text{ and } \quad 
D_2=\dfrac{9g^2}{32\lambda_*^2}+\dfrac{3g^2-\lambda_*^2}{32\lambda_*^2}, $$
satisfies the conditions \begin{equation}
\label{eq:ob4_ir}
 \mathcal{H}[\tilde{h}](q,p) = \mathcal{O}(b^4) \text{ and } 
\mathcal{B}_1[\tilde{h}](q) = 0.\;
 \end{equation}
For the Robin boundary condition we take  
\begin{equation}
\mathcal{B}_0[\tilde{h}](q) = -B_0 \cos q \ b^3 +\mathcal{O}(b^4),
\end{equation}
with the constant $B_0$ given explicitly by the expression 
$$B_0=\frac{1}{2\lambda_*}\dfrac{1}{\left( g^2-\lambda_*^2\right)^{3/2} }\left[ 
\dfrac{\lambda_* p_0 \left( 3g^2-\lambda_*^2\right)^{2}}{g p_0+ 
\lambda_*^{3/2}} 
- \dfrac{g\left( 3g^2+\lambda_*^2\right)^{2}}{2\lambda_*}\right]   .$$

Thus, we propose the following formulation
\begin{equation}
w(q,p)=\tilde{w}(q,p)- B\left(p-p_0\right) \cos q ,
\end{equation}
for some constant $B$. Then the condition \eqref{eq:ob4} changes accordingly. In particular, 
\eqref{eq:ob4_ir}  becomes
 \begin{equation}
 \mathcal{H}[\tilde{h}](q,p) = B \dfrac{p-p_0}{\lambda_*} \cos q \ b^3 + 
\mathcal{O}(b^4)
 \end{equation}
and the Robin boundary condition takes the following form:
\begin{equation}
\mathcal{B}_0[\tilde{h}](q) = \left[ \frac{2}{\lambda_*}\left( g p_0 
+\lambda_*^{3/2}\right) B-B_0 \right] \cos q \ b^3 +\mathcal{O}(b^4).
\end{equation}
Therefore we can choose some appropriate value of $B$ in order to minimize the 
errors observed in the above conditions. For example, \begin{itemize}
\item $\displaystyle B=0$ ensures $\mathcal{H}[\tilde{h}]=\mathcal{O}(b^4).$
\item $\displaystyle B = B_0 \frac{\lambda_*}{2}\frac{1}{g p_0 
+\lambda_*^{3/2}}$ ensures $\mathcal{B}_0[\tilde{h}]=\mathcal{O}(b^4).$
\end{itemize}

We observe the magnitude of the error of the $b^3$ term, for the operators 
$\mathcal{H}$ and $\mathcal{B}_0$, to be 
\begin{equation}
\mathcal{E}_d := B \dfrac{-p_0}{\lambda_*} \geq 0
\end{equation}
and
\begin{equation}
\mathcal{E}_c := B_0 - \frac{2}{\lambda_*}\left( g p_0 +\lambda_*^{3/2}\right) 
B 
\geq 0,
\end{equation}
respectively, with 
$$\displaystyle 0 \leq B \leq B_0 
\frac{\lambda_*}{2}\frac{1}{g p_0 +\lambda_*^{3/2}} .$$

\subsection{Constant vorticity}

While the situation is similar for this case, it is, however, much more complicated. Indeed, 
we apply the expansion of the form \eqref{sol_perturb_3}, with $H$, $m$ and $u$ 
given by \eqref{H*0_r}, \eqref{m0_r} and \eqref{u_r}, respectively, to the 
system \eqref{eq:ob4} and we get the following BVP on the rectangular $R$, for 
the even and periodic function $w(q,p)$:
\begin{equation}\label{eq:w_r}
\begin{aligned}
w_{pp}(q,p)+\dfrac{1}{r^2(p)} w_{qq}(q,p) -\frac{3\gamma}{r^2(p)}w_p(q,p)&= 
f_1(p) \cos q + f_3(p) \cos(3q) , \\
 w_{p}(q,0) - \frac{g}{\lambda_*^{3/2}}w(q,0)  &=  g_1\cos q +g_3 \cos(3q), \\
w(q,p_0)&=0,
\end{aligned}
\end{equation}
where \begin{itemize}
\item $f_1(p)$ and $f_3(p)$ are known functions explicitly dependent on $r(p)$,
\item $g_1$ and $g_3$ are known constants explicitly dependent on the 
parameters 
$\gamma,\ \lambda_*,\ p_0$ and $g$,
\end{itemize}
for which we avoid to write the complete expressions for matters of brevity. 
Proceeding as in Section 3.2, we reduce the solution of the above BVP to the 
following formulation: Determine $w(q,p)$, of the form
\begin{equation} \label{w_b3}
w(q,p)=\sum_{n=0}^\infty w_n(p) \cos(nq),
\end{equation}
with $w_n(p)$ given by
\begin{equation} \label{transform_4}
w_n(p)=\dfrac{\hat{w}_n(r(p))}{r(p)} \qquad \text{ with } \ r(p)=\sqrt{\lambda_*-2\gamma p}\,,
\end{equation}
and where $\hat{w}_n(r)$ are the solutions of the BVPs
\begin{equation}\label{eq:w_r_n-hat}
\begin{aligned}
\dfrac{d^2}{dr^2}\hat{w}_n(r)-\dfrac{n^2}{\gamma^2} \hat{w}_n(r) &= \hat{f}_n(r) , \qquad &r\in\left(r_0,\sqrt{\lambda_*}\right)\\
 \dfrac{d}{dr}\hat{w}_n + \dfrac{g - \gamma \sqrt{\lambda_*}}{\gamma\lambda_*}\hat{w}_n  &=  \hat{g}_n, \qquad &r=\sqrt{\lambda_*}\\
\hat{w}_n&=0, \qquad &r=r_0=\sqrt{\lambda_*-2\gamma p_0},
\end{aligned}
\end{equation}
for some known functions $\hat{f}_n(r)$ and constants $\hat{g}_n$. 

\begin{remark} The functions $\hat{f}_n(r)$ are bilinear combinations of 
hyperbolic functions and polynomials in $1/r$. The constants 
$\hat{g}_n$ are dependent on $\gamma$, $g$, $\lambda_*$ and $r_0$. \end{remark}

Eventually, using the arguments of the previous subsection we derive the 
following formula for $w(q,p)$:
\begin{equation}
w(q,p)=\tilde{w}(q,p)- \tilde{B}r^2(p_0)\dfrac{H(p)}{r(p)}\cos q ,
\end{equation}
for some constant $\tilde{B}$ and $\tilde{w}(q,p)$ given by
\begin{equation}\label{w-r}
\begin{aligned}
\tilde{w}(q,p)=&\left[  \tilde{A}_1(p) \sinh\big( H(p) \big) +\tilde{A}_2(p) 
\sinh\big( 3H(p) \big)\right. \\ 
 & \left. +\tilde{B}_1(p) \cosh\big( H(p) \big) + \tilde{B}_2(p) \cosh\big( 
3H(p) \big)\right] \cos q \\
+&\left[  \tilde{A}_3(p) \sinh\big( H(p) \big) +\tilde{A}_4(p) \sinh\big( 3H(p) 
\big)\right. \\ 
 & \left. +\tilde{B}_3(p) \cosh\big( H(p) \big) + \tilde{B}_4(p) \cosh\big( 
3H(p) \big)\right] \cos(3q) ,
\end{aligned}
\end{equation}
where $\{\tilde{A}_j(p),\tilde{B}_j(p)\}_{j=1}^4$ are functions of the form 
$$\dfrac{1}{r^5(p)}\sum_{n=0}^4 a_{j,n} r^n(p) \quad \text{ and } \quad 
\dfrac{1}{r^5(p)}\sum_{n=0}^4 b_{j,n} r^n(p),$$
respectively. Here
$$\begin{array}{l}
a_{1,0}=- \displaystyle\frac{9}{32}\,\gamma ^2 r_0^3, \qquad a_ {1, 1} = 0\,,\\ [0.2cm]
a_ {1, 2} = -\Big\{r_0^3 \left(3 g \Big(g-\gamma  
\sqrt{\lambda _*}-\lambda _*\right) \sqrt{\lambda _*} \left(g-\gamma  
\sqrt{\lambda _*}+\lambda _*\right) \,,\\[0.2cm]
\qquad\qquad +r_0 \Big(-3 g^3-15 g^2 \gamma  
\sqrt{\lambda _*}-5 \gamma  \left(\gamma ^2-\lambda _*\right) \lambda _*^{3/2} +3 
g \lambda _* \left(5 \gamma ^2+\lambda _*\right)\Big)\Big)\Big\}\Big/ \\[0.2cm]
\qquad\qquad\Big\{32 \left(-r_0 
\left(g-\gamma  \sqrt{\lambda _*}\right)+g \sqrt{\lambda _*}\right) \left(g^2-2 
g \gamma  \sqrt{\lambda _*}+\left(\gamma ^2-\lambda _*\right) \lambda 
_*\right)\Big\} \,,\\[0.2cm]
a_ {1, 3} =-\displaystyle\frac{r_0^3 \left(3 g^2-3 g \gamma  \sqrt{\lambda 
_*}+\left(\gamma ^2-3 \lambda _*\right) \lambda _*\right)}{32 \lambda _*^{5/2}}, 
\qquad a_ {1, 4} =0\,,\\[0.2cm]
b_ {1, 0}=0, \qquad b_ {1, 1} = - \displaystyle\frac{9 \gamma  r_0^3}{32}, \qquad b_ {1, 2} 
= -\frac{\gamma  r_0^3 \left(3 g^2-3 g \gamma  \sqrt{\lambda _*}+\left(\gamma 
^2-3 \lambda _*\right) \lambda _*\right)}{32 \lambda _*^{5/2}} \,,\\[0.2cm]
b_ {1, 3} = \displaystyle\frac{r_0^4 \sqrt{\lambda _*} \left(-3 g^2+3 g \gamma  
\sqrt{\lambda _*}+\lambda _* \left(-\gamma ^2+\lambda _*\right)\right)}{4 
\left(g r_0-\left(g+\gamma  r_0\right) \sqrt{\lambda _*}\right) \left(g^2-2 g 
\gamma  \sqrt{\lambda _*}+\left(\gamma ^2-\lambda _*\right) \lambda _*\right)}\,,\\[0.2cm]
b_ {1, 4} = \displaystyle\frac{r_0^3 \sqrt{\lambda _*} \left(3 g^2-3 g \gamma  
\sqrt{\lambda _*}+\left(\gamma ^2-\lambda _*\right) \lambda _*\right)}{4 \left(g 
r_0-\left(g+\gamma  r_0\right) \sqrt{\lambda _*}\right) \left(g-\gamma  
\sqrt{\lambda _*}-\lambda _*\right) \left(g-\gamma  \sqrt{\lambda _*}+\lambda 
_*\right)}\,,
\end{array}$$

$$\begin{array}{l}
a_{2,0}=\displaystyle\frac{3}{32} \gamma ^2 r_0^3,\qquad a_{2,1}=0,\qquad 
a_{2,2}=\frac{9 r_0^3}{32}\,,\\[0.33cm]
a_{2,3}=\displaystyle\frac{3 r_0^3 \left(3 g^2-3 g \gamma  
\sqrt{\lambda _*}+\left(\gamma ^2-3 \lambda _*\right) \lambda _*\right)}{32 
\lambda _*^{5/2}}, \qquad a_ {2, 4} =0\,,\\[0.33cm]
b_ {2, 0} =0, \qquad b_ {2, 1} =\displaystyle\frac{9}{32} \gamma  r_0^3\,,\\[0.33cm] 
b_ {2, 2} 
=\displaystyle\frac{\gamma  r_0^3 \left(3 g^2-3 g \gamma  \sqrt{\lambda _*}+\left(\gamma 
^2-3 \lambda _*\right) \lambda _*\right)}{32 \lambda _*^{5/2}}, \qquad b_ {2, 3} 
=b_ {2, 4} = 0\,,\\[0.33cm]
a_{3,0}= -\displaystyle\frac{3}{32} \gamma ^2 r_0^3, \quad a_{3,1}=0, \quad 
a_{3,2}=-\frac{r_0^3}{32}\,,\\[0.33cm] 
a_{3,3}= -\displaystyle\frac{r_0^3 \left(3 g^2-3 g \gamma  
\sqrt{\lambda _*}+\left(\gamma ^2-3 \lambda _*\right) \lambda _*\right)}{32 
\lambda _*^{5/2}}, \quad a_{3,4}=0\,,\\[0.33cm]
b_ {3, 0} =0, \qquad b_{3,1}= -\displaystyle\frac{3 \gamma  r_0^3}{32}\,,\\[0.33cm] 
b_{3,2}=-\displaystyle\frac{\gamma  r_0^3 \left(3 g^2-3 g \gamma  \sqrt{\lambda 
_*}+\left(\gamma ^2-3 \lambda _*\right) \lambda _*\right)}{32 \lambda _*^{5/2}}, 
\qquad b_ {3, 3} =b_ {3, 4} = 0\,,\\[0.33cm]
a_{4,0}=\displaystyle\frac{1}{32} \gamma ^2 r_0^3,\qquad a_{4,1}=0,\qquad a_{4,2}=\frac{3 r_0^3}{32}\,,\\[0.33cm] 
a_{4,3}=\displaystyle\frac{3 r_0^3 \left(3 g^2-3 g \gamma  \sqrt{\lambda 
_*}+\left(\gamma ^2-3 \lambda _*\right) \lambda _*\right)}{32 \lambda 
_*^{5/2}}\,,\\[0.33cm]
a_{4,4}= \Big\{r_0^3 \Big(9 g^5-27 g^4 \gamma  
\sqrt{\lambda _*}+11 g^3 \left(3 \gamma ^2-2 \lambda _*\right) \lambda _*-g^2 
\gamma  \left(21 \gamma ^2-50 \lambda _*\right) \lambda _*^{3/2} \\
\qquad\qquad +g \lambda _*^2 
\left(7 \gamma ^4-35 \gamma ^2 \lambda _*+13 \lambda _*^2\right)-\gamma  \lambda 
_*^{5/2} \left(\gamma ^4-9 \gamma ^2 \lambda _*+15 \lambda 
_*^2\right)\Big)\Big\} \Big/\\[0.33cm]
\qquad\qquad \Big\{64 \left(g-\gamma  \sqrt{\lambda _*}\right) \lambda 
_*^5\Big\} \,,\\[0.33cm]
b_ {4, 0} =0, \qquad b_{4,1}=\displaystyle\frac{3}{32} \gamma  r_0^3\,,\\[0.33cm]
b_{4,2}=\displaystyle\frac{\gamma  r_0^3 \left(3 g^2-3 g \gamma  \sqrt{\lambda 
_*}+\left(\gamma ^2-3 \lambda _*\right) \lambda _*\right)}{32 \lambda _*^{5/2}},
\qquad b_ {4, 3} =b_ {4, 4} = 0\,,
\end{array}$$
where we used the abbreviation $r(p_0)=r_0$. Consequently, the conditions \eqref{eq:ob4} 
take on the following form:
\begin{equation}
 \mathcal{H}[\tilde{h}](q,p) = \tilde{B} r^2(p_0)\dfrac{H(p)}{r^3(p)} \cos q \ 
b^3 + \mathcal{O}(b^4)
 \end{equation}
and 
\begin{equation}
\mathcal{B}_0[\tilde{h}](q) = \left[ 2r^2(p_0)\frac{2 g p_0 + 
r(p_0)\sqrt{\lambda_*}\left( r(p_0)+\sqrt{\lambda_*}\right)}{\lambda_*^{3/2} 
\left( r(p_0)+\sqrt{\lambda_*}\right) } \tilde{B}-\tilde{B}_0 \right] \cos q \ 
b^3 +\mathcal{O}(b^4),
\end{equation}
with 
$$\tilde{B}_0 = D \dfrac{\sum_{n=0}^5 c_n \gamma^n}{ g^2-\lambda_*^2 
-\gamma\sqrt{\lambda_*}\left( 2g-\gamma\sqrt{\lambda_*}\right) },$$
where $$D=-\frac{1}{4\lambda_*^{7/2}}
\dfrac{r(p_0)+\sqrt{\lambda_*}}{2 g p_0 + r(p_0)\sqrt{\lambda_*}\left( 
r(p_0)+\sqrt{\lambda_*}\right) }$$
and 
$$\begin{array}{l}
c_0  = 3g \sqrt{\lambda_*}\left( 9g^4-14g^2 
\lambda_*^2+5\lambda_*^4\right)-4g\sqrt{\lambda_*}r_0\left( 
9g^4-9g^2\lambda_*^2+4\lambda_*^4\right)\\[0.33cm]
\qquad\qquad -\displaystyle\dfrac{2p_0}{\sqrt{\lambda_*}+r_0} \left( 9g^6-12g^4 \lambda_*^2 + 13 g^2 
\lambda_*^4 - 2\lambda_*^6\right)\,,\\[0.33cm]
c_1=-\lambda _*^{3/2} \left(33 g^4-39 g^2 \lambda _*^2+4 \lambda 
_*^4\right)-r_0 \left(-60 g^4 \lambda _*+45 g^2 \lambda _*^3-7 \lambda 
_*^5\right)\,,\\[0.2cm]
c_2=-3 \left(-7 g^3 \lambda _*^2+5 g \lambda _*^4+2 g r_0 \lambda _*^{3/2} 
\left(9 g^2-4 \lambda _*^2\right)\right)\,,\\[0.2cm]
c_3=-\lambda _*^{5/2} \left(7 g^2-2 \lambda _*^2\right)-r_0 \left(-28 g^2 
\lambda _*^2+5 \lambda _*^4\right)\,,\\[0.2cm]
c_4=-8 g r_0 \lambda _*^{5/2}+g \lambda _*^3, \qquad c_5=r_0 \lambda _*^3\,.
\end{array}$$
Similarly to the previous section we observe that the magnitude of the error of 
the $b^3$ term, for the operators $\mathcal{H}$ and $\mathcal{B}_0$, to be 
\begin{equation}
\mathcal{E}_d := \tilde{B} \dfrac{H(p_0)}{r(p_0)} \geq 0
\end{equation}
and
\begin{equation}
\mathcal{E}_c := \tilde{B}_0 - 2r^2(p_0)\frac{2 g p_0 + 
r(p_0)\sqrt{\lambda_*}\left( r(p_0)+\sqrt{\lambda_*}\right)}{\lambda_*^{3/2} 
\left( r(p_0)+\sqrt{\lambda_*}\right) } \tilde{B}\geq 0,
\end{equation}
respectively, with 
$$\displaystyle 0 \leq \tilde{B} \leq \tilde{B}_0 
\frac{\lambda_*^{3/2}}{2r^2(p_0)}\frac{r(p_0)+\sqrt{\lambda_*}}{2 g p_0 + 
r(p_0)\sqrt{\lambda_*}\left( r(p_0)+\sqrt{\lambda_*}\right)} .$$

\begin{remark} One can see that the case $\gamma=0$, equivalent to 
$r(p)=\sqrt{\lambda_*}$ for all $p\in[p_0,0]$, in the above expressions yields the 
formulae of the previous subsection. \end{remark}

\section{Illustration of the solutions}
Considering the figures depicted in this section, we recall that the value of the gravitational 
constant of acceleration is $g=9.8$ and we fix the relative mass flux $p_0=-2$. 
This section is organised as follows. In the first subsection we illustrate the 
results for the second order asymptotic expansion. In the second subsection we 
present the analogue results for the third order expansions. In the third 
subsection we present some figures that compare the obtained second order approximation 
results against existing results i.e. first 
order asymptotic expansions. Finally, in the fourth subsection we compare the 
results between the second and the third order approximations.

The illustrations below are given for several values of the vorticity $\gamma$, 
as indicated in each figure. Moreover, each streamline is defined by $[p=\hbox{constant}]$, for 
constants ranging from $p_0$ to $0$. We note that for the wave height 
$h$ and the velocity field $v$, we have that: 
\begin{itemize}
\item $p=p_0$ gives $h=0$ and $v=0$.
\item $h$ is an increasing function with respect to $p$.
\item $v$ is an increasing function with respect to $p$, for $q>0$, and a 
decreasing function with respect to $p$, for $q<0$.
\end{itemize}
whereas, for the pressure $P$:
\begin{itemize}
\item $P(q,0)=P_{atm}(q)=P_{atm}$. 
\item $P$ is a decreasing function with respect to $p$.
\end{itemize}

Let us recall the definition of the wave height as the maximal oscillation of the free boundary, 
$h(0,0)-h(\pi,0)$, with the understanding the the wave crest is located at $q=0$. Moreover, we recall 
from \cite{CS} that:
 \begin{itemize}
\item The velocity field is given by  $\displaystyle \quad (c-u,v)=\left( \frac{1}{h_p},-\frac{h_q}{h_p}\right).$
\item The water pressure is given by 
$$\quad P=P_{atm}-\frac{1+h_q^2}{2h_p^2}- g h - \gamma p+\dfrac{Q}{2},$$
for some  constant $Q$.
\end{itemize}

\subsection{Second order}
In this subsection we depict the results that we obtain from the second order 
asymptotic expansion of the solution $h(q,p)$, given in 
\eqref{sol_perturb_2_r}. 
Concerning the other parameters, we remind that $\lambda_*$ varies with the 
vorticity according to the dispersion relation \eqref{disp}.

Firstly, {\textit{for each value of the constant vorticity $\gamma$,}} 
we have to fix the value of $b$ by allowing a small error on condition 
\eqref{eq:ob3_r}, i.e.
$$\max_b\left\{\left\Vert\mathcal{H}\right\Vert_2,\left\Vert\mathcal{B}
_0\right\Vert_2\right\}=\epsilon \text{ for some small } \epsilon>0.$$

\begin{figure}[ht!]
\centering
\includegraphics[scale=1.3]{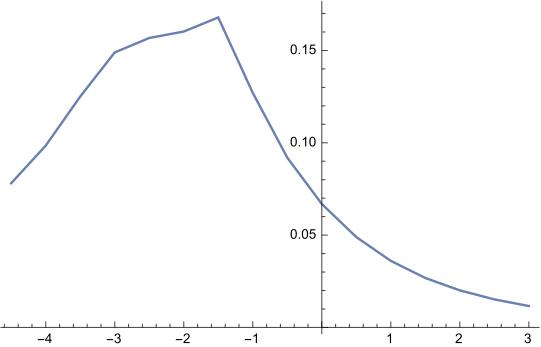}
\caption{The value of $b$ for different values of vorticity, varying from 
$\gamma=-4.5$ to $\gamma=3$.}
\label{fig:b-tb0}
\end{figure}

In Figure \ref{fig:b-tb0} we show how $b$ varies for different values of 
vorticity. The value $b$ is increasing when $\gamma\in[-4.5, -1.5]$, with  
$\left\Vert\mathcal{H}\right\Vert_2$ dominant in the above definition. On the contrary, 
the value of $b$ is decreasing when $\gamma\in[-1.5, 3]$, with $\left\Vert\mathcal{B}_0\right\Vert_2$ dominant in 
the above definition.

\vfill\eject

\subsubsection{Wave profiles}
\begin{figure}[ht!]
\centering
\includegraphics[scale=0.85]{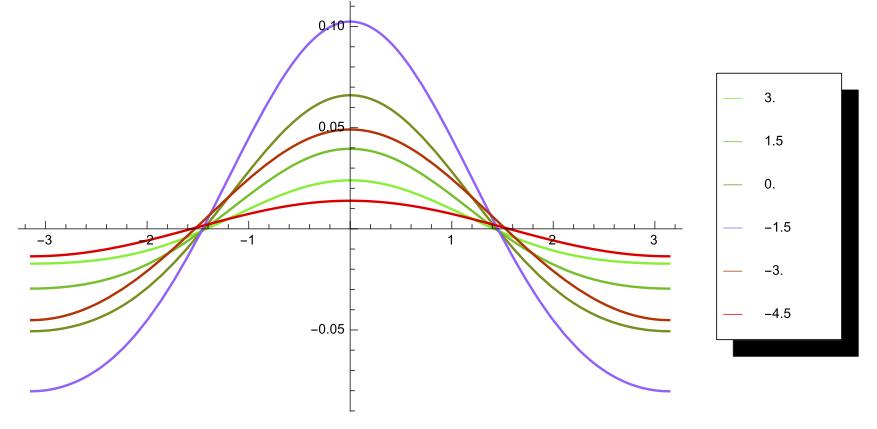}
\caption{The free surface for different values of vorticity, varying from 
$\gamma=-4.5$ to $\gamma=3$. Here the wave profile $\eta(x)$ is illustrated, which is  
oscillating around $y=0$.}
\label{fig:free2}
\end{figure}

\begin{figure}[ht!]
\centering
\includegraphics[scale=1]{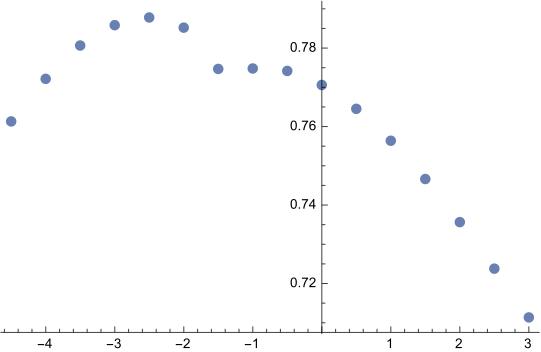}
\caption{The mean depth for different values of vorticity, varying from $\gamma=-4.5$ to $\gamma=3$.}
\label{fig:d2}
\end{figure}

\begin{figure}[ht!]
\centering
\includegraphics[scale=1]{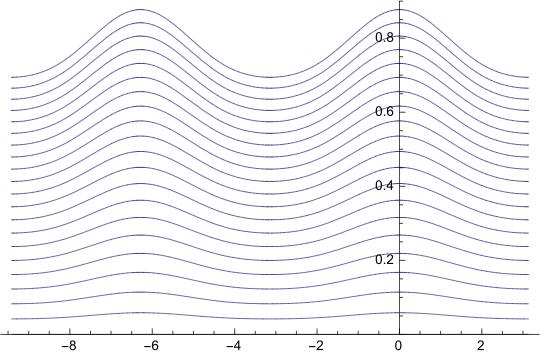}
\caption{The height of the water along streamlines for $\gamma=-1.5$, over two 
wavelengths.}
\label{fig:str2}
\end{figure}

\begin{figure}[ht!]
\centering
\includegraphics[scale=1]{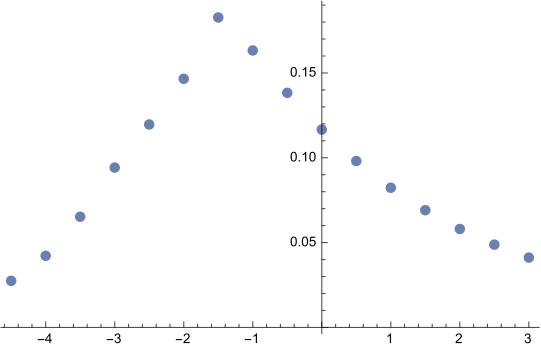}
\caption{The wave height for different values of vorticity, 
varying from $\gamma=-4.5$ to $\gamma=3$.}
\label{fig:amp2}
\end{figure}

\subsubsection{Velocities}

\begin{figure}[ht!]
\centering
\includegraphics[scale=1.25]{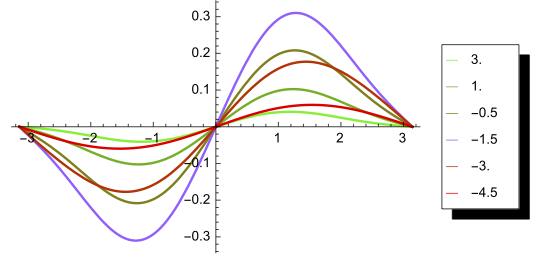}
\caption{The vertical velocity $v$ on free surface for different values of 
vorticity, varying from $\gamma=-4.5$ to $\gamma=3$.}
\label{fig:vfree2}
\end{figure}

\begin{figure}[ht!]
\centering
\includegraphics[scale=1.25]{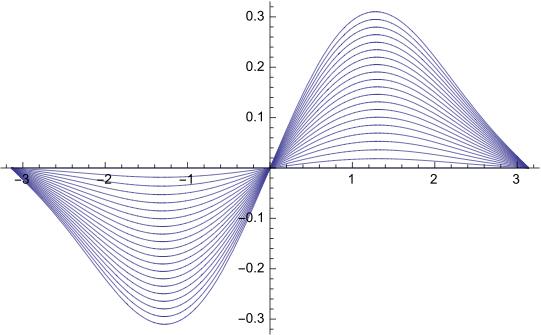}
\caption{The vertical velocity $v$ of the flow along streamlines for 
$\gamma=-1.5$.}
\label{fig:v2}
\end{figure}

\begin{figure}[ht!]
\centering
\includegraphics[scale=1]{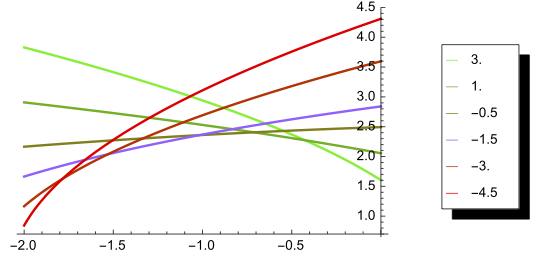}
\caption{The horizontal velocity $c-u$ on the straight line below the crest, i.e. 
$q=0$ and $p\in [p_0,0]$, for different values of vorticity, varying from 
$\gamma=-4.5$ to $\gamma=3$.}
\label{fig:u2}
\end{figure}

\subsubsection{Pressure}

\begin{figure}[ht!]
\centering
\includegraphics[scale=1.33]{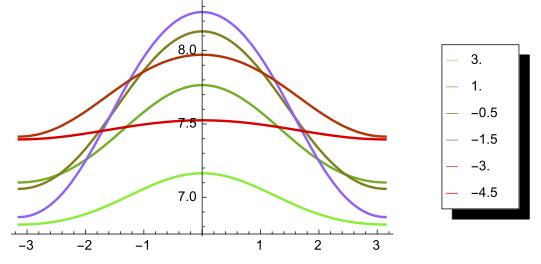}
\caption{The water pressure on the bottom, $p=p_0$, for different values of 
vorticity, varying from $\gamma=-4.5$ to $\gamma=3$.}
\label{fig:prbot2}
\end{figure}

\begin{figure}[ht!]
\centering
\includegraphics[scale=1.25]{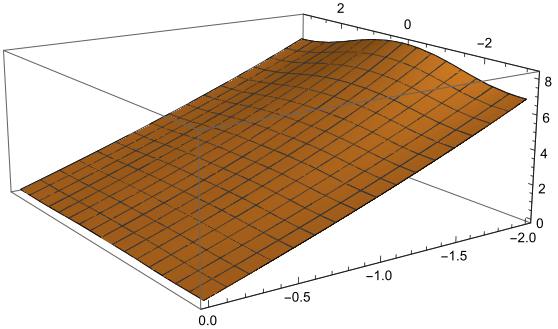}
\caption{The water pressure of the flow for $\gamma=-1.5$: ($P-P_{atm}$) 
vanishes on the free boundary $p=0$ and takes the maximal 
value on the bottom $p=p_0$.}
\label{fig:pr2}
\end{figure}

\subsection{Third order}
Here, we proceed as follows: $\left\Vert\mathcal{H}\right\Vert_2$ and 
$\left\Vert\mathcal{B}_0\right\Vert_2$ for the third order asymptotic expansion 
are, in general, functions of $\gamma$, $b$ and $\tilde{B}$, see \eqref{sol_perturb_3} and \eqref{w-r}. So, for each value 
of $\gamma$, we pick the specific value of $b$ that we have obtained in the 
previous subsection,  as illustrated in Figure \ref{fig:b-tb0}. Then, 
for each pair $(\gamma,b)$, we determine the value of $\tilde{B}$ that minimizes the quantity
$$\left\Vert\mathcal{H}\right\Vert_2^2+\left\Vert\mathcal{B}_0\right\Vert_2^2.$$

\begin{figure}[ht!]
\centering
\includegraphics[scale=0.7]{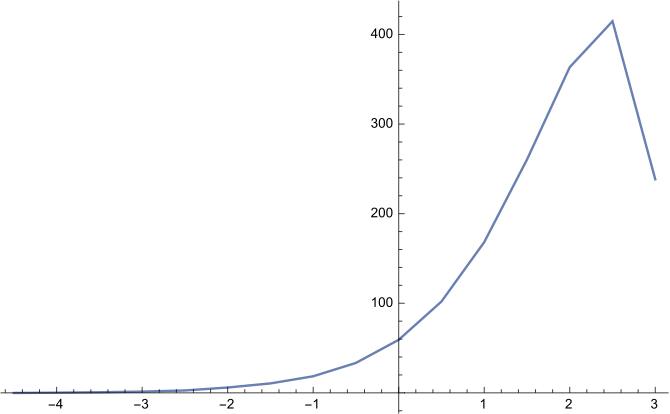}
\caption{The value of $\tilde{B}$ for different values of vorticity, varying 
from $\gamma=-4.5$ to $\gamma=3$.}
\label{fig:a-ta0}
\end{figure}

In Figure \ref{fig:a-ta0} we show how $\tilde{B}$ varies for different values of 
vorticity. In what follows we will display the analogue of the 
figures of the previous subsection, but now for the third order asymptotic 
expansion \eqref{sol_perturb_3}.

\subsubsection{Wave profiles}


\begin{figure}[ht!]
\centering
\includegraphics[scale=1]{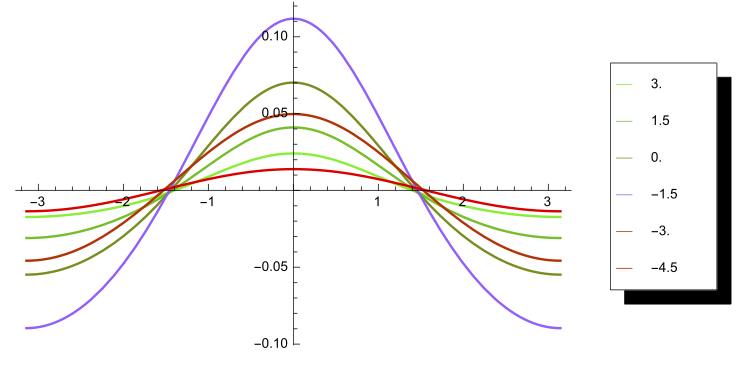}
\caption{The free boundary for different values of vorticity, varying from 
$\gamma=-4.5$ to $\gamma=3$.}
\label{fig:free3}
\end{figure}

\begin{figure}[ht!]
\centering
\includegraphics[scale=1]{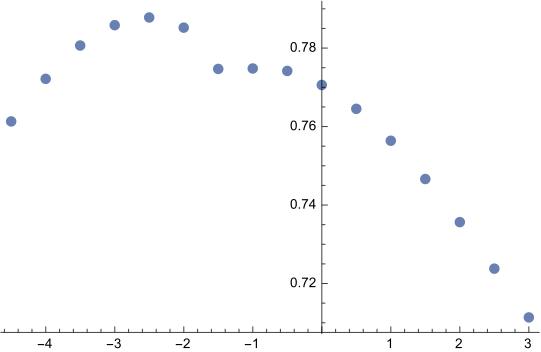}
\caption{The mean depth for different values of vorticity, varying from  $\gamma=-4.5$ to $\gamma=3$.}
\label{fig:d3}
\end{figure}

\begin{figure}[ht!]
\centering
\includegraphics[scale=1]{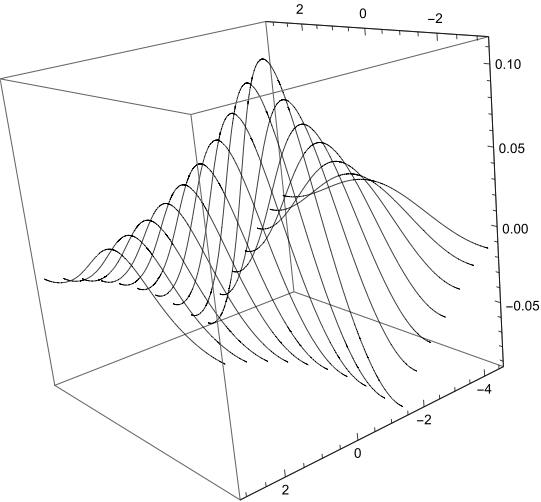}
\caption{Each of the disjoint curves displays the wave profile $\eta(x)$, for different 
values of vorticity, varying from $\gamma=-4.5$ to $\gamma=3$.}
\label{fig:free3-3d}
\end{figure}

\begin{figure}[ht!]
\centering
\includegraphics[scale=1]{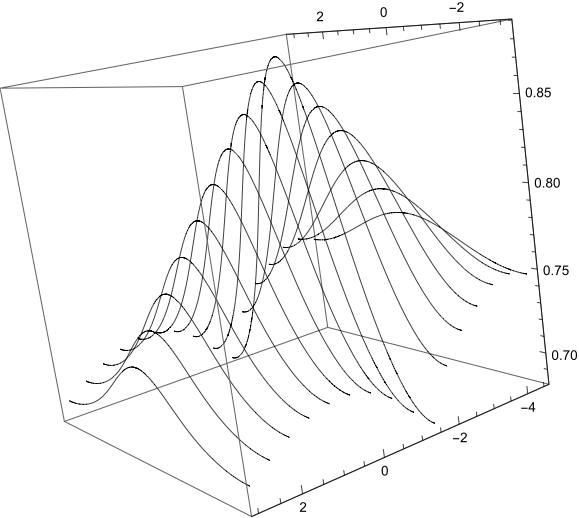}
\caption{Each of the disjoint curves displays the free boundary $h(q,0)=\eta(x)+d$, for different 
values of vorticity, varying from $\gamma=-4.5$ to $\gamma=3$. This illustration, compared to Fig. \eqref{fig:free3-3d}, takes into account the mean depth $d$, which varies with the vorticity $\gamma$, 
see Fig. \ref{fig:d3}.}
\label{fig:free3d}
\end{figure}

\begin{figure}[ht!]
\centering
\includegraphics[scale=1]{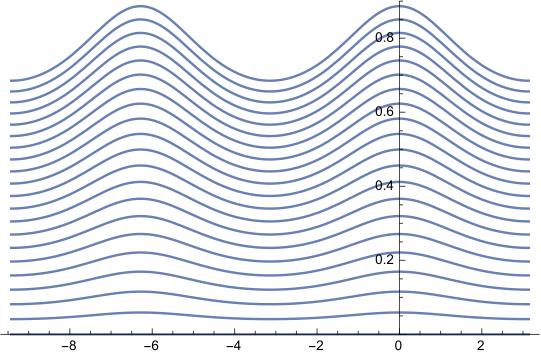}
\caption{The height of the water along streamlines for $\gamma=-1.5$, over two 
wavelengths.}
\label{fig:str3}
\end{figure}

\begin{figure}[ht!]
\centering
\includegraphics[scale=1]{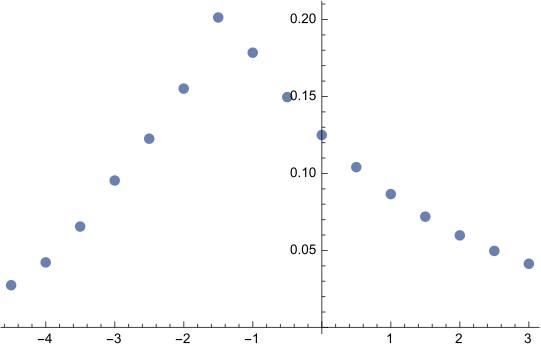}
\caption{The wave height for different values of the vorticity, 
varying from $\gamma=-4.5$ to $\gamma=3$.}
\label{fig:amp3}
\end{figure}

\subsubsection{Velocities}

\begin{figure}[ht!]
\centering
\includegraphics[scale=1.25]{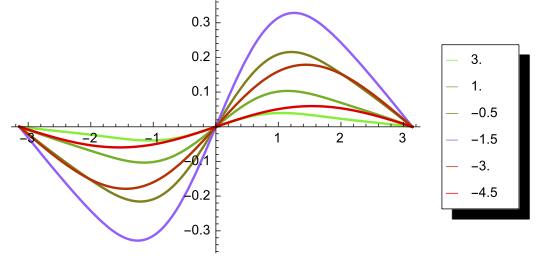}
\caption{The vertical velocity $v$ on free surface for different values of 
vorticity, varying from $\gamma=-4.5$ to $\gamma=3$.}
\label{fig:vfree3}
\end{figure}

\begin{figure}[ht!]
\centering
\includegraphics[scale=1.25]{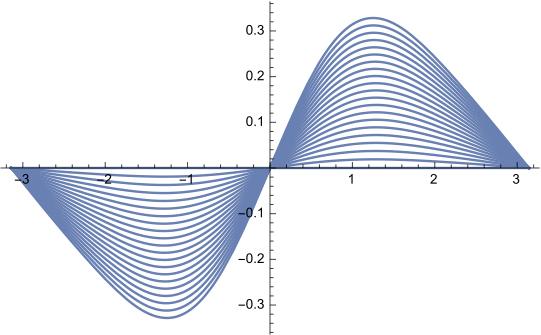}
\caption{The vertical velocity $v$ of the flow along streamlines for 
$\gamma=-1.5$.}
\label{fig:v3}
\end{figure}

\begin{figure}[ht!]
\centering
\includegraphics[scale=1]{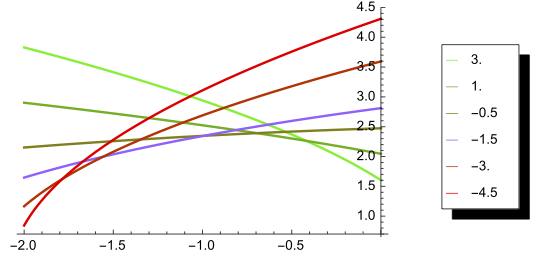}
\caption{The horizontal velocity $c-u$ on the straight line below the crest, i.e. 
$q=0$ and $p \in [p_0,0]$, for different values of vorticity, varying from 
$\gamma=-4.5$ to $\gamma=3$.}
\label{fig:u3}
\end{figure}

\subsubsection{Pressure}

\begin{figure}[ht!]
\centering
\includegraphics[scale=1]{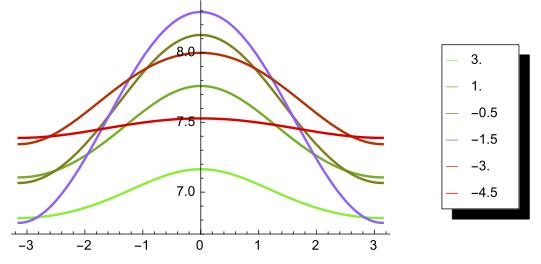}
\caption{The water pressure on the bottom, $p=p_0$, for different values of 
vorticity, varying from $\gamma=-4.5$ to $\gamma=3$.}
\label{fig:prbot3}
\end{figure}


\begin{figure}[ht!]
\centering
\includegraphics[scale=0.7]{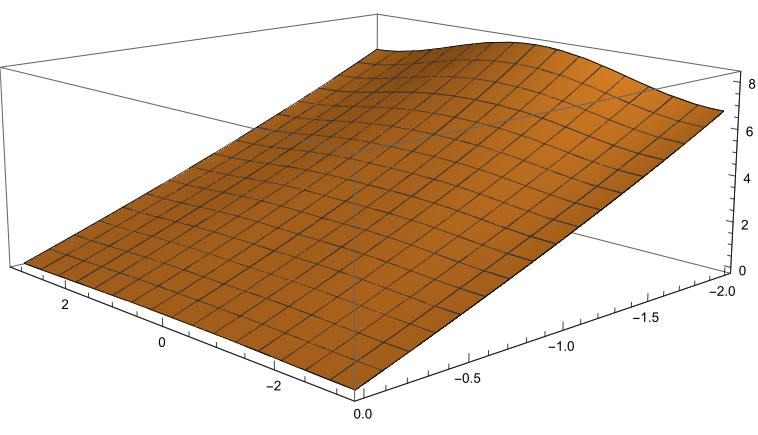}
\caption{The water pressure of the flow for $\gamma=-1.5$: ($P-P_{atm}$) 
vanishes on the free boundary $p=0$ and takes the maximal 
value on the bottom $p=p_0$.}
\label{fig:pr3}
\end{figure}

\subsection{Comparison with existing results}

\begin{figure}[ht!]
\centering
\includegraphics[scale=1.25]{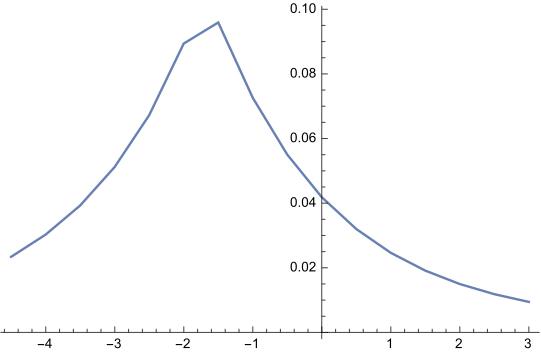}
\caption{The value of $b$ for different values of vorticity, varying from 
$\gamma=-4.5$ to $\gamma=3$.}
\label{fig:b-tb1}
\end{figure}

Here we compare the results that we obtained from the second order asymptotic 
expansion with the ones already existing in the literature, i.e. the first 
order 
asymptotic expansion. In order to choose the value of $b$ for the depiction of 
the first order asymptotic expansion, we proceed as in Section 5.1 and we obtain 
qualitative results similar to Figure \ref{fig:b-tb0}, but now $b$ is in general smaller, 
as depicted in Figure \ref{fig:b-tb1}.

\begin{figure}[ht!]
\centering
\includegraphics[scale=1]{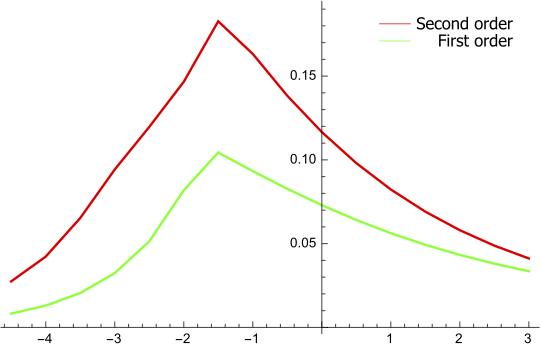}
\caption{The height of wave for different values of vorticity, 
varying from $\gamma=-4.5$ to $\gamma=3$.}
\label{fig:amp12}
\end{figure}

This choice of $b$ indicates the comparison on the wave height for different values of vorticity, depicted in Figure \ref{fig:amp12}.  There we observe that for large values of the absolute value of vorticity, $|\gamma|$, the wave height is relatively small and the first order approximation shows similar results with the second order one. However, the shape of the water profile differs significantly, see Fig. \ref{fig:str12-p} and \ref{fig:curv23}. Moreover, for the rest values of vorticity, we observe a discrepancy in the wave height between the two approximations. 

\subsubsection{Negative vorticity $\gamma=-1.5$}

In what follows we fix our vorticity to $\gamma=-1.5$ in order to illustrate the 
improvement provided by higher-order asymptotics.

\begin{figure}[ht!]
\centering
\includegraphics[scale=1]{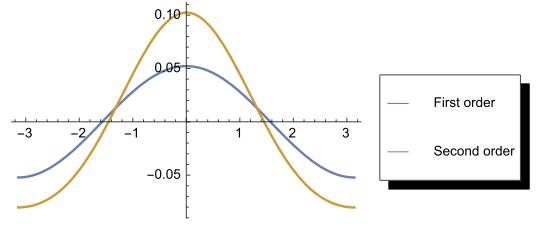}
\caption{The free surface for vorticity  $\gamma=-1.5$. The mean depth is 0.80072 and 0.79222 for the first and the second order approximation, respectively.}
\label{fig:free12-m}
\end{figure}

\begin{figure}[ht!]
\centering
\begin{subfigure}{.49\textwidth}
\centering
\includegraphics[scale=0.67]{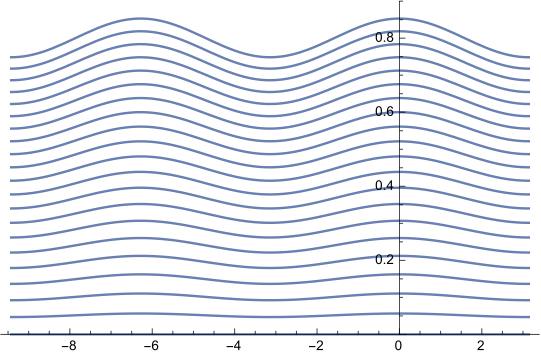}
\subcaption{First order asymptotic.}
\label{fig:str1-12-m}
\end{subfigure}
\begin{subfigure}[ht!]{.49\textwidth}
\centering
\includegraphics[scale=0.67]{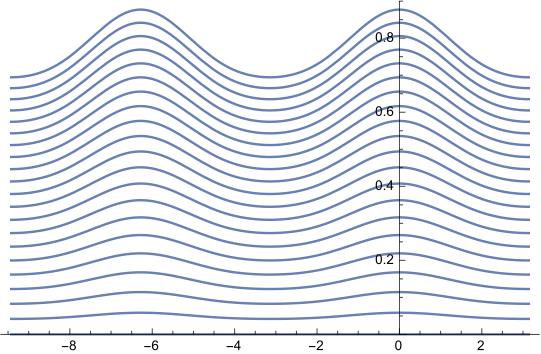}
\subcaption{Second order asymptotic.}
\label{fig:str2-12-m}
\end{subfigure}
\caption{The height of the water along streamlines for $\gamma=-1.5$, over two 
wavelengths.}
\label{fig:str12-m}
\end{figure}

\begin{figure}[ht!]
\centering
\includegraphics[scale=1]{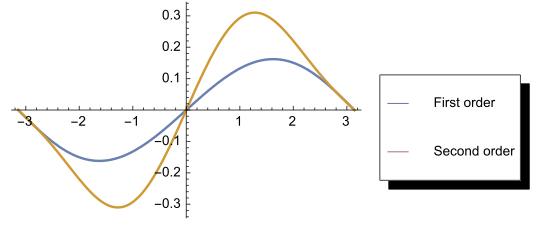}
\caption{The vertical velocity $v$ on the free surface, for vorticity  
$\gamma=-1.5$.}
\label{fig:vfree12-m}
\end{figure}

\begin{figure}[ht!]
\centering
\includegraphics[scale=1]{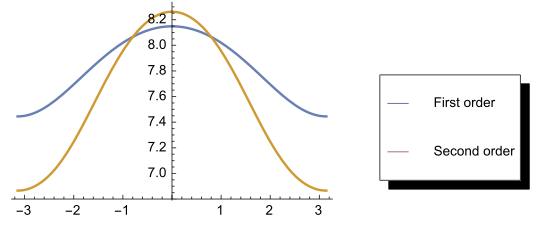}
\caption{The water pressure: $P-P_{atm}$ on the bottom, for the choice of vorticity  
$\gamma=-1.5$.}
\label{fig:prbot12-m}
\end{figure}

\subsubsection{Positive vorticity $\gamma=1.5$}

In what follows we fix our vorticity to $\gamma=1.5$ in order to illustrate the 
effect of higher-order asymptotics.

\begin{figure}[ht!]
\centering
\includegraphics[scale=1]{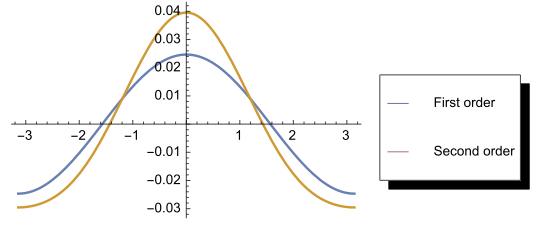}
\caption{The free surface for vorticity  $\gamma=1.5$. The mean depth is 0.758042 and 0.752228 for the first and the second order approximation, respectively.}
\label{fig:free12-p}
\end{figure}

\begin{figure}[ht!]
\centering
\begin{subfigure}{.49\textwidth}
\centering
\includegraphics[scale=0.65]{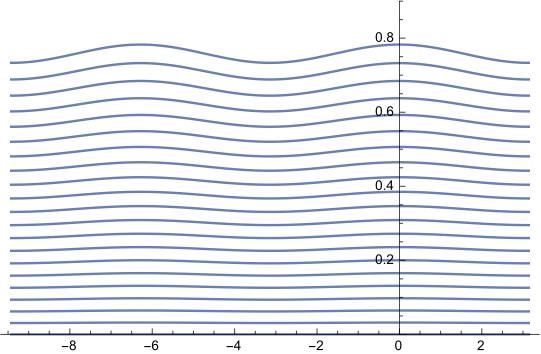}
\subcaption{First order asymptotic.}
\label{fig:str1-12-p}
\end{subfigure}
\begin{subfigure}[ht!]{.49\textwidth}
\centering
\includegraphics[scale=0.65]{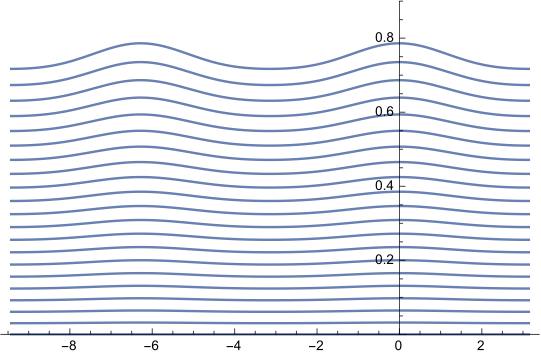}
\subcaption{Second order asymptotic.}
\label{fig:str2-12-p}
\end{subfigure}
\caption{The height of the water along streamlines for $\gamma=1.5$, over two 
wavelengths.}
\label{fig:str12-p}
\end{figure}

\begin{figure}[ht!]
\centering
\includegraphics[scale=1]{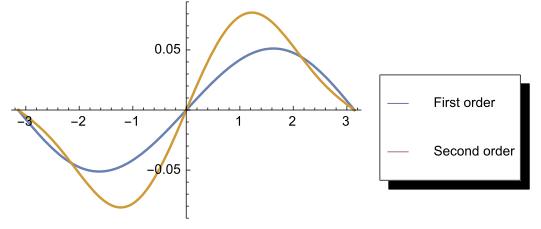}
\caption{The vertical velocity $v$ on the free surface, for vorticity  
$\gamma=1.5$.}
\label{fig:vfree12-p}
\end{figure}

\begin{figure}[ht!]
\centering
\includegraphics[scale=1]{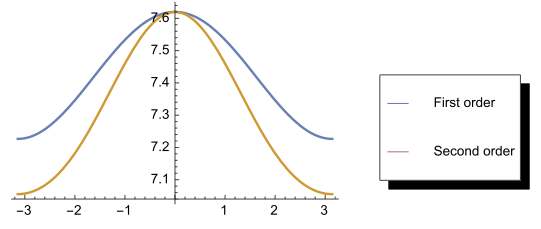}
\caption{The water pressure: $(P-P_{atm})$ on the bottom, for vorticity  
$\gamma=1.5$.}
\label{fig:prbot12-p}
\end{figure}

\subsection{Comparison of second and third order asymptotics}

In what follows we fix our vorticity to $\gamma=1.5$ and we compare the 
difference between the second and third order expansion.

\begin{figure}[ht!]
\centering
\includegraphics[scale=1]{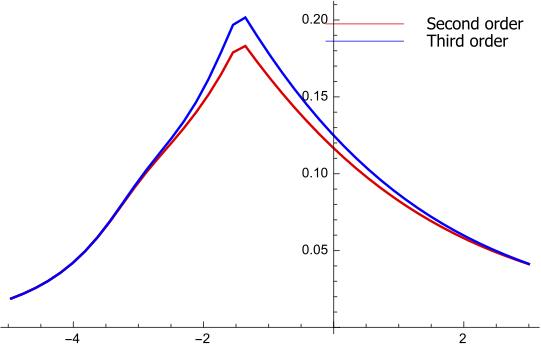}
\caption{The wave height for different values of vorticity, 
varying from $\gamma=-4.5$ to $\gamma=3$.}
\label{fig:amp23}
\end{figure}

\begin{figure}[ht!]
\centering
\begin{subfigure}{.49\textwidth}
\centering
\includegraphics[scale=0.6]{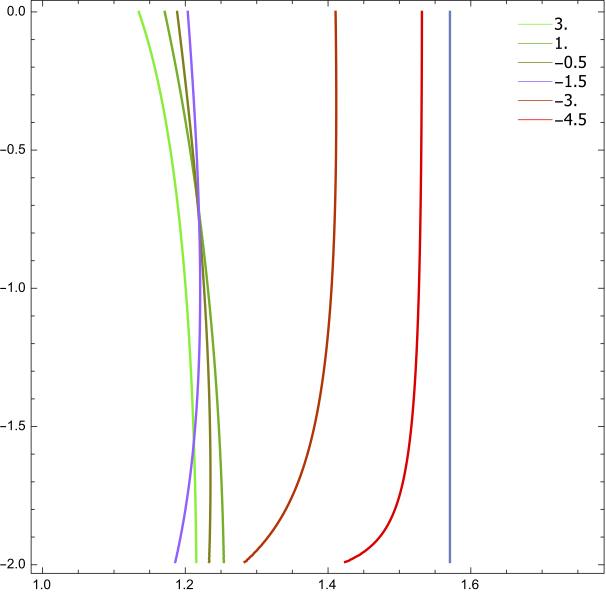}
\subcaption{Second order asymptotic.}
\label{fig:curv2}
\end{subfigure}
\begin{subfigure}[ht!]{0.49\textwidth}
\centering
\includegraphics[scale=0.6]{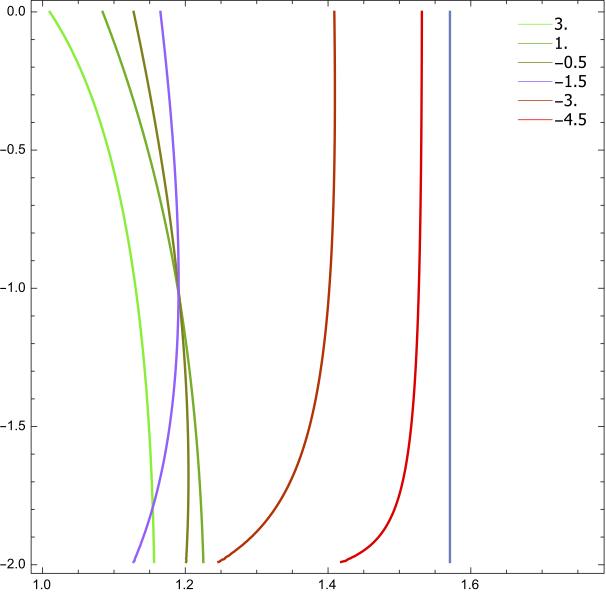}
\subcaption{Third order asymptotic.}
\label{fig:curv3}
\end{subfigure}
\caption{The curve tracking the change of curvature, for values of the vorticity 
varying from $\gamma=-4.5$ to 
$\gamma=3$.}
\label{fig:curv23}
\end{figure}

In Figure \ref{fig:curv23} we depict, for different values of vorticity, the curves $h_{qq}=0$ for 
$p\in[p_0,0]$, thus indicating where the curvature changes. Indeed, let us fix the value of vorticity $\gamma$. 
If for instance $p=\tilde{p}$ is fixed, then the equation $$\dfrac{\partial^2h(q;\tilde{p})}{\partial q^2}=0$$ has a unique solution $q=\tilde{q},$ where $ \tilde{q}\in(0,\pi/2)$. Allowing $p$ varying from $p_0$ to $0$, we get the curve $p\mapsto \tilde{q}(p)$. For different values of the vorticity we get the different curves depicted in Figure \ref{fig:curv23}, observing where the wave profile changes from convex to concave.

\section{Conclusions}

Among the conclusions that we can extract from our asymptotic expansions and the associated 
figures, the most important are the following:
\begin{itemize} 
\item We observe a considerable improvement from the first to the second asymptotic expansion, noticeable in 
the depictions of the wave height in Fig. \ref{fig:amp12} as well as in the depictions of the profile of the wave in 
Fig. \ref{fig:str12-m} and \ref{fig:str12-p}.
\item We observe a small difference in wave height from the second to the third asymptotic expansion, cf. Fig. \ref{fig:amp23}. 
However, there is a considerable difference in the change of the curvature, see Fig. \ref{fig:curv23}. This is due to the choice of the parameters $b$ and $\tilde{B}$. In the particular examples we did choose to keep the same values for $b$ and to minimize the errors $\Vert\mathcal{H}\Vert_2$ and 
$\Vert\mathcal{B}_0\Vert_2$, meaning that we kept approximately the same wave height and 
improved the accuracy of the wave profile.
\item Fig. \ref{fig:free3} shows that an opposing current has a steepening effect on the wave profile. As already pointed 
out in the Introduction, first-order linear theory does not capture this important feature of wave-current 
interactions. 
\end{itemize}

Asymptotic expansions beyond those obtained here will be presented in upcoming publications.

\section*{Acknowledgements}
K. Kalimeris was supported by the project {\it Computation of large amplitude water waves} (P 27755-N25), funded by the Austrian Science Fund (FWF).

\end{document}